\documentclass[12pt]{amsart}
\usepackage{amscd,amssymb}
\usepackage[all]{xy}

\title{Motivic integration over Deligne-Mumford stacks}
\author{Takehiko Yasuda}
\address{Graduate School of Mathematical Sciences, University of Tokyo,
Komaba, Meguro, Tokyo, 153-8914, Japan}

\email{t-yasuda@ms.u-tokyo.ac.jp}

\theoremstyle{plain}
\newtheorem{thm}{Theorem}[section]

\newtheorem{prop}[thm]{Proposition}
\newtheorem{cor}[thm]{Corollary}
\newtheorem{lem}[thm]{Lemma}

\theoremstyle{definition}
\newtheorem{defn}[thm]{Definition}

\theoremstyle{remark}
\newtheorem{rem}[thm]{Remark}

\def\AA{\mathbb A}
\newcommand{\CC}{\mathbb C}

\newcommand{\QQ}{\mathbb Q}

\newcommand{\ZZ}{\mathbb Z}
\newcommand{\LL}{\mathbb L}

\newcommand{\G}{\mathbf{G}}

\newcommand{\NN}{\mathbb N}

\newcommand{\cA}{\mathcal{A}}

\newcommand{\cC}{\mathcal{C}}
\newcommand{\cD}{\mathcal{D}}
\newcommand{\cE}{\mathcal{E}}

\newcommand{\cG}{\mathcal{G}}
\newcommand{\cI}{\mathcal{I}}
\newcommand{\cJ}{\mathcal{J}}
\newcommand{\cK}{\mathcal{K}}

\newcommand{\cO}{\mathcal{O}}

\newcommand{\cU}{\mathcal{U}}
\newcommand{\cV}{\mathcal{V}}
\newcommand{\cX}{\mathcal{X}}
\newcommand{\cY}{\mathcal{Y}}
\newcommand{\cZ}{\mathcal{Z}}
\newcommand{\cW}{\mathcal{W}}

\newcommand{\fm}{\mathfrak{m}}
\newcommand{\fn}{\mathfrak{n}}
\newcommand{\fs}{\mathfrak{s}}
\newcommand{\fI}{\mathfrak{I}}

\newcommand{\fR}{\mathfrak{R}}
\newcommand{\fS}{\mathfrak{S}}

\newcommand{\Isom}{\mathcal I som}

\newcommand{\Ker}{\mathrm{Ker}\,}
\newcommand{\Coker}{\mathrm{Coker}\,}

\newcommand{\ord}{\mathrm{ord}\,}

\newcommand{\Hom}{\mathrm{Hom}}
\newcommand{\Spec}{\mathrm{Spec}\,}

\newcommand{\GL}{\mathrm{GL}}
\newcommand{\SL}{\mathrm{SL}}
\def\Im{\mathrm{Im}}
\newcommand{\Conj}{\mathrm{Conj}}

\newcommand{\Gr}{\mathrm{Gr}}
\newcommand{\Aut}{\mathrm{Aut}}

\newcommand{\rank}{\mathrm{rank}\,}

\newcommand{\id}{\mathrm{id}}
\newcommand{\codim}{\mathrm{codim}\,}
\newcommand{\diag}{\mathrm{diag}}

\newcommand{\sing}{\mathrm{sing}}
\newcommand{\chara}{\mathrm{char}\,}

\newcommand{\Jac}{\mathrm{Jac}}
\newcommand{\sht}{\mathrm{sht}}

\newcommand{\Gal}{\mathrm{Gal}}

\newcommand{\discrep}{\mathrm{discrep}}
\newcommand{\age}{\mathrm{age}}

\newcommand{\bmu}{\boldsymbol{\mu}}
\newcommand{\Affk}{(\mathrm{Aff}/k)}

\numberwithin{equation}{section}

\begin{document}

\maketitle

\begin{abstract}
The aim of this article is to develop the theory of motivic integration over Deligne-Mumford stacks and 
to apply it to the birational geometry of Deligne-Mumford stacks. 
\end{abstract}

\tableofcontents

\section{Introduction}
In this article, we study the motivic integration over Deligne-Mumford stacks,
which was started in \cite{Twistedjet}.
The motivic integration was introduced by Kontsevich \cite{Orsay} and
developed by Denef and Loeser \cite{germs}, \cite{DL-quotient} etc. 
It is now well-known that the motivic integration is
 effective in the study of birational geometry.
  For example, Batyrev \cite{non-archi} has applied it to the study of stringy $E$-functions and
  Musta\c{t}\v{a} \cite{mustata-jet-inventiones}
 to one of the singularities appearing in the minimal model program. 

We first recall the motivic integration over varieties. 
Thanks to Sebag \cite{Sebag}, we can work over an arbitrary perfect field $k$. 
Let $X$ be a variety over $k$, that is, a separated algebraic space
of finite type over $k$. For a non-negative integer $n$, 
an {\em $n$-jet} of $X$ over  a $k$-algebra $R$ is a $R[[t]]/t^{n+1}$-point of $X$.  
For each $n$, there exists an algebraic space $J_n X$  parameterizing $n$-jets. 
For example, $J_0 X$ is $X$ itself and $J_1 X$ is the tangent bundle of $X$.
The spaces $J_n X$, $n \in  \ZZ_{\ge 0}$ constitute a projective system and the limit
$ \displaystyle J_\infty X := \lim _{\longleftarrow} J_n X$ exists.
We can define a measure $\mu_X$ and construct an integration theory on $J_\infty X$ with values in some ring (or semiring) in which we can add and multiply 
the classes $\{V\}$ of varieties $V$ and some class of infinite sums are defined.
For example, we can use a completion of the Grothendieck ring of mixed Hodge structures
($k=\CC$) or mixed Galois representations ($k$ a finite field).
If $X$ is smooth, then we have
\[
 \int _{J_\infty X} 1 \ d\mu_X = \mu_X (J_\infty X) = \{X\} .
\]

To generalize the theory to Deligne-Mumford stacks, it is not sufficient to consider 
only $R[[t]]/t^{n+1}$-points  of a stack.
Inspired by a work of Abramovich and Vistoli \cite{AV},
the author introduced the notion of  twisted jets in \cite{Twistedjet}.
Let $\cX$ be a separated Deligne-Mumford stack of finite type over $k$ and $\bmu_{l,k}$
be the group scheme of 
$l$-th roots of unity for a positive integer $l$ prime to the characteristic of $k$. 
A {\em twisted $n$-jet} over $\cX$ is a representable morphism
from a quotient stack $[ (\Spec R[[t]]/t^{nl +1}) /\bmu_{l,k}] $ to $\cX$.
We will prove that the category $\cJ_n \cX$  of twisted $n$-jets 
is a Deligne-Mumford stack.
If $k$ is algebraically closed and
 $\cX$ is a quotient stack $[M/G]$, then we have
\[
 \cJ_0 \cX \cong \coprod _{g \in \Conj (G)} [M^g /C_g] .
\] 
Here $\Conj (G)$ is a representative set of conjugacy classes, $M^g$ the fixed point locus of $g$ and
$C_g$ is the centralizer of $g$. The right hand side often appears in the study of McKay correspondence.
There exists also the projective limit $\displaystyle \cJ_\infty \cX := \lim _{\longleftarrow} \cJ_n \cX$. 
When $\cX$ is smooth, we define a measure $\mu_\cX$ and construct an integration theory on the point set $|\cJ_\infty \cX|$.

Let $\LL$ be the class $\{\AA^1_k\}$ of an affine line.
To a variety $X$ and an ideal sheaf $\cI \subset \cO _X$, we can associate a function
 $\ord \cI : J_\infty X \to \ZZ_{\ge 0} \cup \{\infty\}$
 and a function $\LL^{\ord \cI}$.
Consider a proper birational morphism $f : Y \to X$ of varieties with $Y$ smooth.
The Jacobian ideal sheaf $\Jac _f \subset  \cO_Y$ is defined to be the $0$-th Fitting ideal of $\Omega _{Y/X}$. If $X$ is also smooth, then this is identical
with the ideal sheaf of the relative canonical divisor $K_{Y/X}:= K_Y - f^* K_X$.
Let $f_\infty : J_\infty Y \to J_\infty X$ be the morphism induced by $f$.
 The relation of the measures $\mu_X$ and $\mu_Y$ is described by the following {\em transformation rule}:
\[
\int F d\mu_X = \int  (F\circ f_\infty) \LL ^{ - \ord \Jac_f}d\mu_Y .
\]
 This formula was proved by Kontsevich  \cite{Orsay},  Denef and Loeser \cite{germs}, and Sebag \cite{Sebag}. 
Using this, we obtain many results in the birational geometry. 
For instance,
Kontsevich proved the following: 
If $f: Y \to X$ and $f' :Y \to X'$ are proper birational morphisms of smooth proper varieties over $\CC$,
and if $ K_{Y/X} = K_{Y/X'} $, 
then 
the Hodge structure of $H^i (X ,\QQ)$ and that of $H^i (X',\QQ)$ are isomorphic.

We generalize the transformation rule to Deligne-Mumford stacks.
If we consider only representable morphisms, no interesting phenomenon appears.
A morphism of Deligne-Mumford stacks is said to be {\em birational} if it induces an isomorphism of open dense substacks.
 For example, if $M$ is a variety with an effective action of a finite group $G$, then the natural morphism from the quotient stack $[M/G]$ to the quotient variety $M/G$ is
  birational. 
A morphism $f:\cY \to \cX$ is said to be {\em tame} if for every geometric point $y $ of $\cY$,
 $\Ker (\Aut (y) \to \Aut (f(x)))$ is of order prime to the 
characteristic of $k$.
The transformation rule is generalized to tame, proper and birational morphisms.

Let $\tilde x$ be a geometric point of $\cJ_0 \cX$ and $x$ its image in $\cX$.
A $\bmu_l$-action on the tangent space $T_x \cX$ derives from $\tilde x$.
If for  suitable basis, $\zeta \in \bmu_l$ acts by 
$\diag (\zeta ^ {a_1} ,\dots , \zeta^{ a_d})$, $1 \le a_i \le l$, then we define 
\[
\sht (\tilde x) := d - \frac{1}{l} \sum _{i=1}^l a_i .
\]
Thus we have a function $\sht : |\cJ_0 \cX| \to \QQ$. 
We denote by $\fs_\cX$, the composite $|\cJ_\infty \cX| \to |\cJ_0 \cX| \xrightarrow{\sht} \QQ$.

Also for a birational morphism $f:\cY \to \cX$ of Deligne-Mumford stacks,
we define its Jacobian ideal sheaf, $\Jac_f$, in the same way.
However, the associated order function $\ord \Jac _f$ is a $\QQ$-valued function.

Now we can formulate the generalized transformation rule as follows:

\begin{thm}
Let $f:\cY \to \cX$ be a tame proper birational morphism of Deligne-Mumford stacks
of finite type and pure dimension.
Suppose that $\cY$ is smooth and that $\cX$ is either a smooth Deligne-Mumford stack or
a variety. Then we have
\[
\int F \LL^{ \fs _\cX } d\mu_{\cX} = 
\int (F \circ f_\infty) \LL ^{ - \ord \Jac _f + \fs _\cY} d\mu_{\cY} .
\]
(See  Theorem \ref{thm-transformationstackstack} for  details).
\end{thm}

\begin{rem}
The theorem was proved in \cite{Twistedjet} for the morphisms from a smooth stack $\cX$ without reflection
to its coarse moduli space with Gorenstein singularities. 
\end{rem}

We apply the transformation rule to the birational
geometry of Deligne-Mumford stacks. 
We recall Batyrev's
work \cite{non-archi} as a background of this subject. 
Suppose $k=\CC$. 
Let $M$ be a smooth variety, $G$ a finite group acting effectively on $M$
and $X = M/G$ be the quotient variety. 
By calculations, 
Batyrev showed a relation of the orbifold $E$-function of the $G$-variety $M$
and the stringy $E$-function of $X$.
Denef and Loeser \cite{DL-quotient} proved a similar result with motivic integration. 
From the viewpoint of stacks, the orbifold $E$-function is defined rather for the quotient stack
$[M/G]$ than for the $G$-variety $M$. 
The natural morphism $[M/G] \to X$ is proper and  birational. 
Then Batyrev's result can be viewed as a relation of invariants of birational stacks.
We will reformulate his results in a full generality from this viewpoint.
 
Let $\cX$ be a smooth Deligne-Mumford stack of finite type over a perfect field $k$, 
$D = \sum u_i D_i $ a $\QQ$-divisor on $\cX$
and $W \subset |\cX|$ a constructible subset. 
We put $\fI_D := \sum u_i \cdot \ord \cI_{D_i}$.
We define an invariant
\[
\Sigma _W(\cX,D):=\int _{\pi^{-1}(W)} \LL^{ \fI_D  + \fs_\cX} d\mu_\cX .
\]
Here $\pi:\cJ_\infty \cX \to \cX $ is the natural projection.
The function $\sht:|\cJ_0 \cX| \to \QQ$ is, in fact, locally constant and for a connected component
$\cV$, $\sht (\cV) \in \QQ$ is well-defined.
If $D=0$ and $W =|\cX| $, then the invariant is equal to
\[
 \sum_{\cV \subset \cJ_0 \cX} \{ \cV\}\LL^{\sht (\cV)} .
\]
If $k=\CC$ and $\cX$ is proper, this has the information of the Hodge structure of 
the orbifold cohomology defined below.

In characteristic zero, we can generalize the invariant $\Sigma _W (\cX,D)$ 
to the case where $\cX$ is singular:
A {\em log Deligne-Mumford stack} is defined to be the pair $(\cX ,D)$
of a normal Deligne-Mumford stack $\cX$ of finite type and a $\QQ$-divisor $D$ on $\cX$ such that
 $K_\cX +D$ is $\QQ$-Cartier.
For a log Deligne-Mumford stack $(\cX ,D)$
and a constructible subset $W\subset |\cX|$,
if $f:\cY \to \cX$ is a proper birational morphism with $\cY$ smooth,
then we define $\Sigma _W (\cX , D) := \Sigma_{f^{-1}(W)} (\cY , f^* (K_\cX +D)-K_\cY)$.
This invariant is a common generalization and refinement of the stringy $E$-function and the orbifold $E$-function.
By a calculation, we will see that 
 $\Sigma_W (\cX ,D) \ne \infty$ if and only if $(\cX ,D)$ is
Kawamata log terminal around $W$ (For the definition, see Definition \ref{def-KLT}).

The following is the direct consequence of the transformation rule and viewed as a generalization 
of Batyrev's result and Denef and Loeser's one.

\begin{thm}
Let $(\cX,D)$ and $(\cX',D')$ be log Deligne-Mumford stacks. 
Assume that there exist a smooth DM stack $\cY$ and proper birational morphisms $f:\cY \to \cX$ and 
$f': \cY \to \cX$ such that $f^* (K_\cX +D ) = (f')^* (K_{\cX'}+D')$
 and $f^{-1}(W) = (f')^{-1}(W')$.
 In positive characteristic,  assume in addition that $\cX$ and $\cX'$ are smooth and that
 $f$ and $f'$ are tame.
Then we have 
\[
 \Sigma _W(\cX ,D)= \Sigma _{W'}(\cX',D') .
\]
\end{thm}

\begin{rem}
Kawamata \cite{kawamata} obtained a closely related result in terms of the derived category.
\end{rem}

Finally we give corollaries of this theorem.

Let $G \subset \GL _d( \CC) $ be a finite subgroup and $X := \CC^d/G$ the quotient variety.
For $g \in G$, we define a rational number $\age (g)$ as follows: 
Let $l$ be the order of $g$ and $\zeta := \exp (2 \pi \sqrt{-1}/l)$.
If we write $g = \diag (\zeta ^{a_i} ,\dots , \zeta ^{a_d}) $, $ 0 \le a_i \le l-1$,
for suitable basis of $\CC$, then 
\[
 \age (g ) := \frac{1}{l} \sum _{i=1} ^d a_i .
\] 
If $g \in \SL_d(\CC)$ , then $\age (g)$ is an integer.
The following was called the Homological McKay correspondence.
It was  proved by Y.\ Ito and Reid \cite{ito-reid}
for dimension three and by Batyrev for arbitrary dimension \cite{non-archi}.
(See also \cite{reid-bourbaki}).

\begin{cor}
Suppose that $G \subset \SL_d(\CC)$ and that 
there is a  crepant resolution $Y \to X$.
For an even integer $i$, put 
\[
n_i := \sharp \{ g \in \Conj (G) | \age (g) = i/2 \}.
\]
Then we have
\[
 H^i (Y , \QQ) \cong 
 \begin{cases}
0  & (i: \text {odd})\\
\QQ(-i/2)^{ \oplus n_i} & (i:\text{even}) .
\end{cases}
\]
\end{cor}

Since $X = \CC^d /G$ has only quotient singularities, $K_X$ is $\QQ$-Cartier and its pull-back by arbitrary morphism
is defined.
For a resolution $f :Y \to X $ and for each exceptional prime divisor $E \subset Y$, 
there is a rational number $a(E,X)$ such that
\[
 K_Y \equiv f^* K_ X + \sum_{E \subset Y} a (E,X) E .
\] 
The {\em discrepancy} of $X$ is defined to be the infimum of $a(E ,X)$
for all resolutions $Y \to X$ and all exceptional divisors $E \subset Y$. 
The following is a refinement of Reid--Shepherd-Barron--Tai  criterion
for canonical (or terminal) quotient singularities (see \cite[\S 4.11]{reid}).

\begin{cor}
For a finite group  $G \subset \GL_d (\CC) $ without reflection, the discrepancy of $X = \CC^d /G$ is
equal to 
\[
 \min \{ \age (g) | 1 \ne g \in G\} -1.
\]
\end{cor}

Chen and Ruan \cite{CR} defined a new cohomology for topological orbifolds (Satake's $V$-manifolds),
called {\em orbifold cohomology}. We give its algebraic version.
Let $\cX$ be a smooth Deligne-Mumford stack over $\CC$.
For $i \in \QQ$, we define 
\[
 H^i _{orb} (\cX , \QQ) := \bigoplus _{\cV \subset \cJ_0 \cX} H ^ {i-2 \sht (\cV)} (\bar \cV , \QQ) \otimes \QQ(- \sht (\cV)) .  
\]
Here $\bar \cV$ is the coarse moduli space of $\cV$. If $\cX$ is proper, then
$H^i_{orb} (\cX , \QQ)  $ is  a pure Hodge structure of weight $i$. 
(We define Hodge structure with fractional weights in the trivial fashion.)
The following was conjectured by Ruan \cite{ruan} and a weak version was proved by Lupercio-Poddar \cite{Lupercio-Poddar} and the author \cite{Twistedjet} independently.
This is a generalization of Kontsevich's theorem stated above.

\begin{cor}
Let $\cX $ and $\cX'$ be proper and smooth Deligne-Mumford stacks of finite type over $\CC$.
Suppose that there exist a smooth Deligne-Mumford stack $\cY$ and
 proper birational morphisms $f : \cY \to \cX $ and $ f' : \cY \to \cX'$
such that $K_{\cY /\cX} = K_{\cY /\cX'}$.
Then for every $i \in \QQ$, there is an isomorphism of Hodge structures
\[
 H^i _{orb}(\cX ,\QQ) \cong H^i _{orb }( \cX' , \QQ) . 
\]
\end{cor}

We also define the $p$-adic orbifold cohomology.
Let $\cX$ be a smooth Deligne-Mumford stack over a finite field $k$.
and $p$  a prime number different from the characteristic of $k$.
If necessary, replacing $k$ with its finite extension, 
we define
\[
 H _{orb}^i (\cX \otimes \bar k , \QQ_p) :=\bigoplus _{\cV \subset \cJ_0 \cX} 
 H^{i -2 \sht (\cV)} (\bar \cV \otimes \bar k , \QQ_p) \otimes \QQ_p (-\sht (\cV)) .
\]
Replacing $k$ is necessary to ensure that fractional Tate twists
$\QQ_p (-\sht (\cV))$ exist. 

\begin{cor}
Let $\cX $ and $\cX'$ be proper and smooth Deligne-Mumford stacks of finite type over a finite field $k$.
Suppose that there exist a smooth Deligne-Mumford stack $\cY$ and
 tame proper birational morphisms $f : \cY \to \cX $ and $ f' : \cY \to \cX'$
such that $K_{\cY /\cX} = K_{\cY /\cX'}$.
Suppose that the $p$-adic orbifold cohomology groups of $\cX$ and $\cX'$ are defined.
Then for every $i \in \QQ$, there is an isomorphism of Galois representations
\[
 H^i _{orb}(\cX \otimes _k \bar k,\QQ_p)^{ss} \cong H^i _{orb }( \cX'\otimes _k \bar k , \QQ_p)^{ss} . 
\]
Here the superscript ``$ss$'' means the semisimplification.
\end{cor}

For varieties,  T.\ Ito \cite{ito-tetsu-Hodge} and Wang \cite{wang} obtained a similar result over number fields.

\subsection{Notation and convention}
Throughout this paper, we work over a perfect base field $k$.
A Deligne-Mumford stack (DM stack for short) is supposed to be separated.
What we mean by a variety is a separated algebraic space of finite type over $k$.

\begin{itemize}
\item $\NN$, $\ZZ_{\ge 0}$ : the set of positive integers and that of non-negative integers
\item $[M/G]$ : quotient stack
\item $|\cX|$ : the set of points of $\cX$
\item $\bar \cX$ : the coarse moduli space of a DM stack $\cX$
\item $\cD_{n,S} ^l := [D_{nl,S} /\bmu_{l,k} ]  $ 
\item $D_{n,S} := \Spec R[[t]]/t^{n+1}$ $( S = \Spec R )$
\item $\bmu_l \subset \bar k$  : the group of $l$-th roots of unity 
\item  $\bmu _{l,k} := \Spec k[x]/(x^l-1)$ : the group scheme of $l$-th roots of unity over $k$
\item $\Conj (G)$ : a representative set of conjugacy classes $[g]$ of $g \in G$
\item $\Conj (\bmu_l , G)$ :  a representative set of conjugacy classes  of $\bmu_l \hookrightarrow G$
\item $J_nX$: $n$-jet space
\item $J_n^{(a)} X$ : For a scheme with $G$-action and $a : \bmu_l \hookrightarrow G$,
$J_n^{(a)} X \subset J_n X$ is the locus where the two $\bmu_l$-actions on $J_n X$ coincide
\item $\cJ^l_n \cX$ : the stack of twisted $n$-jets of order $l$ 
\item $\cJ_n \cX := \coprod _{\chara(k) \nmid  l} \cJ^l_n \cX$ : the stack of twisted $n$-jets
\item $\pi_n:\cJ_\infty \cX \to \cJ_n \cX$, $\pi :\cJ_\infty \cX \to \cX$ : natural projections
\item $f_n :\cJ_n \cY \to \cJ_n\cX $ : the morphism induced by $f : \cY \to \cX$
\item $\fR$, $\fS$ : the semirings of equivalence classes of convergent stacks and convergent spaces
\item $\LL:= \{\AA^1_k\}$
\item $MHS$ and $MHS^{1/r}$ : the category of mixed Hodge structures and the category
 of $\frac{1}{r}\ZZ$-indexed ones
\item $\G_k := \Gal (\bar k/k)$ : absolute Galois group
\item $MR( \G_k , \QQ_p)$ : the category of mixed Galois representations
\item $\mu _\cX$ : motivic measure
\item $\sht (p)$, $\sht (\cV)$ : shift number
\item $\fs _\cX$ : the composite $\cJ_\infty \cX \xrightarrow{\pi_0} \cJ_0 \cX \xrightarrow{\sht} \QQ$
\item $\Jac _f$, $\Jac _X$ : the Jacobian ideal sheaves of a morphism $f$ and a variety $X$
\item $\ord \cI$ : the order function of an ideal sheaf $\cI$ over $\cJ_\infty \cX$
\item $\fI _D$ : For a $\QQ$-divisor $D = \sum u_i D_i$, if $\cI_{D_i}$ is the ideal sheaf of $D_i$,
then $\fI_D := \sum u_i \ord \cI_{D_i}$.
\item $\omega _\cX$, $K_\cX$ and $K_{\cY/\cX}$ : 
 canonical sheaf,  canonical divisor and relative canonical divisor
\end{itemize}

\subsection{Acknowledgments}

I am deeply grateful to Yujiro Kawamata and Fran\c{c}ois Loeser for their encouragement and useful advice. 
I would also like to thank  Masao Aoki, Julien Sebag and Orlando Villamayor 
for useful conversations, and
Sindhumathi Revuluri and Laura Sobrin for their proofreading.
The referee's comments were helpful in revising this paper.
I have written this paper during my stay at the \'{E}cole Normale Sup\'{e}rieure.
I was fortunate enough to study in a favorable environment.
Financial support has been provided by the Japan Society for the Promotion of Science.
%%%%%%%%%%%%%%%%%%%%%%%%%%%%%%%%%%%%%%%%%%%%%%%%%%%%%%%%%%%%%%%%%%%%%%%%%%
%%%%%%%%%%%%%%%%%%%%%%%%%%%%%%%%%%%%%%%%%%%%%%%%%%%%%%%%%%%%%%%%%%%%%%%%%%

\section{Stacks of twisted jets}

\subsection{Short review of the Deligne-Mumford stacks}

\subsubsection{}

We first review the Deligne-Mumford (DM) stack very briefly.
We mention the book of Laumon and Moret-Bailly \cite{LMB} as a reference of stacks.
We will sometimes use results from it.

Fix a base field $k$. Let $\Affk$ be the category of affine schemes over $k$.
 A DM stack $\cX$ is a category equipped with a functor $\cX \to \Affk$
which satisfies several conditions. It should be a fibered category over $\Affk$ and is usually best understood in terms of the fiber categories $\cX(S)$, for $S \in \Affk$, and
the pull-back functors $f^*:\cX (T) \to \cX(S)$ for $f:S \to T$. 
The $\cX(S)$ are groupoids with, at least for $S$ of finite type, finite automorphism groups.

The DM stacks constitute a 2-category. 
In terms of the fiber categories, a 1-morphism (or simply morphism) $f:\cY \to \cX$ is the data of functors
$f_S :\cY(S) \to \cX(S)$, compatible with pull-backs, and a 2-morphism $f \to g$ is a system of morphisms
of functors $f_S \to g_S$, compatible with pull-backs.
A scheme, or more generally an algebraic space $X$ is identified with the DM stack with fibers the discrete
categories with sets of objects the $X(S) := \Hom (S,X)$.
A diagram of stacks
\[\xymatrix{
\cX \ar[r] ^f \ar[dr]_h &\cY \ar[d]^g \\
& \cZ 
}\]
is said to be commutative if a  2-isomorphism $g \circ f \cong h$ has been given.
The strict identity $g \circ f = h$
is not required.

A morphism $f:\cY \to \cX$ of DM stacks is called {\em representable} if for every morphism $M \to \cX$ with
$M$ an algebraic space, the fiber product $M \times _\cX \cY$ is also an algebraic space.
It is equivalent to that for every object $\xi \in \cY$, the natural map $\Aut (\xi ) \to \Aut (f(\xi))$
is injective. 
We can generalize many properties of a morphism of schemes to DM stack; \'etale, smooth, proper etc.
By a condition in the definition, for every DM stack $\cX$, there exist an algebraic space $M$ and
 an \'etale surjective morphism $M \to \cX$, which is called an {\em atlas}. We say that $\cX$ is smooth, normal etc if an atlas is so. 

The diagonal morphism $\Delta:\cX \to \cX \times \cX$ of a DM stack $\cX$ is
, by definition, representable. We say that $\cX$ is {\em separated} if $\Delta$
is finite, that is, quasi-finite and proper. 
Note that $\Delta$ is not immersion unless $\cX$ is an algebraic space.
In this paper, every DM stack is supposed to be separated. 

\subsubsection{Points and coarse moduli space}

A {\em point} of a DM stack $\cX$ is an equivalence class of morphisms $\Spec K \to \cX$ with $K \supset k$ a field by
the following equivalence relation;
 morphisms $\Spec K_1 \to \cX$ and $\Spec K_2 \to \cX$ are equivalent if 
 there is another field $K_3 \supset K_1, K_2 \supset k$ making the following diagram commutative.
\[\xymatrix{
\Spec K_3 \ar[r] \ar[d]& \Spec K_2 \ar[d]  \\
  \Spec K_1 \ar[r] & \cX
}
\]
We denote by the set of the points by $|\cX|$. It carries a Zariski topology; 
$A \subset |\cX|$ is an open subset if $A = |\cY|$ for some open immersion $\cY \hookrightarrow \cX$.
(see \cite{LMB} for details). 
If $X$ is a scheme, then $|X|$ is equal to the underlying topological space as sets.

A coarse moduli space of a DM stack $\cX$ is an algebraic space equipped with a morphism $\cX \to X$
such that  every morphism $\cX \to Y$ with $Y$ algebraic space uniquely factors through $X$ and for 
every algebraically closed field $K \supset k$, the map $\cX (K) / isom \to X(K)$ is bijective.
By the definition, it is clear that the coarse moduli space is unique up to isomorphism. 
Keel and Mori \cite{Keel-Mori} proved that for every DM stack, the coarse moduli space exists.
If $X$ is the coarse moduli space of $\cX$, then the map $|\cX|\to |X|$ is a homeomorphism.

\subsubsection{Quotient stack}

One of the simplest examples is the quotient stack. Let $M$ be an algebraic space and $G$ a finite group
(or an \'etale finite group scheme over $k$) acting
on $M$. Then we can define the quotient stack $[M/G]$ as follows; an object over a scheme $S$ is 
a pair of a $G$-torsor $P \to S$ and a $G$-equivariant morphism $P \to M$ and a morphism of 
$(P \to S , P \to M)$ to $(Q \to T , Q \to M)$ over a morphism $S \to T$ is a $G$-equivariant 
morphism $P \to Q$ compatible with the other morphisms.
This stack has the canonical atlas $M \to [M/G]$.
There is also a natural morphism $[M/G] \to M/G$ which makes $M/G$ the coarse moduli space.
The composition $M \to [M/G] \to M/G $ is the quotient map.

\subsection{Stacks of twisted jets}

\subsubsection{}
In the article \cite{Twistedjet}, the author introduced the notion of twisted jets. 
There,  only twisted jets
over fields were considered and  the stack of twisted jets was constructed as a closed substack of another stack.
By contrast, in this paper, we consider the category of twisted jets parameterized by arbitrary affine scheme
and verify that it is actually a DM stack.

We first recall jets and arcs over a variety. Here we mean a separated algebraic space of finite type
by a variety.
Let $X$ be a variety and $n$ a non-negative integer. The functor
\begin{align*}
\Affk &\to (\text{Sets}) \\
 \Spec R & \mapsto \Hom (\Spec R[[t]]/t^{n+1} , X)
\end{align*}
is representable by a variety $J_n X$, called the {\em $n$-jet space}.
The natural surjection $k[[t]]/t^{n+2} \twoheadrightarrow  k[[t]]/t^{n+1}$ induces
a natural projection
$J _{n+1} X \to J_{n} X$. Since they are all affine morphisms,
the projective limit $ J_\infty X := \displaystyle \lim _{\longleftarrow} J_n X $ exists.
This is an algebraic space, but not generally of finite type.
We call this the {\em arc space}.
For every field extension $K \supset k$, there is an identification 
\[
 \Hom (\Spec K, J_\infty X) = \Hom (\Spec K[[t]] , X) .
\]
An {\em arc} of $X$ is a point of $J_\infty X$, that is, a morphism 
$\Spec K[[t]] \to X$.

For $S=\Spec R \in \Affk$ and a non-negative integer $n$, we put
\begin{align*}
D_{n,S} &:= \Spec R[[t]]/t^{n+1}. 
\end{align*}
For $l$ a positive integer prime to the characteristic of $k$, 
we denote by $\bmu_l \subset \bar k$ the cyclic group of $l$-th roots of unity.
We define also the group scheme of $l$-th roots of unity over  $k$
\[
\bmu_{l,k} := \Spec k[x]/(x^l-1) .
\]
When $\bmu_{l,k} $ is a constant group scheme, then we identify it with the group $\bmu_l$.
The natural action of $\bmu_{l,k}$ on $D_{n,S}$  is defined  by $ t \mapsto x \otimes t$.
We put
\[
\cD_{n,S}^l := [D_{nl,S}/ \bmu _{l,k} ].
\] 
Also for $n= \infty$, and for a field $K \supset k$, 
 we put 
 \[
  D_{\infty,K} := \Spec K[[t]] \text{ and } \cD_{\infty , K}^l := [D_{\infty ,K}/ \bmu _{l,k}].
 \]

\begin{defn}
Let $\cX$ be a DM stack.
A {\em twisted $n$-jet of order $l$} of $\cX$ over $S$ is a representable morphism
$\cD_{n,S}^l \to \cX$. For a field $K \supset k$, a {\em twisted arc (or twisted infinite jet) of order $l$}
of $\cX$ over $K$ is a representable morphism $\cD_{\infty, K}^l \to \cX$.
\end{defn}

\begin{defn}
Let $\cX$ be a DM stack. Suppose $n < \infty$.
We define the \textit{stack of twisted $n$-jets of order $l$}, denoted $\cJ _n^l \cX$, as follows;
an object over $S \in \Affk$ is a representable morphism $\cD_{n,S}^l \to \cX$, 
a morphism from $\gamma:\cD_{n,S}^l \to \cX$  to $\gamma':\cD_{n,T}^l \to \cX$ 
over $f:S \to T$ is a 2-morphism from $\gamma$ to $ f' \circ \gamma'$, where $f' : \cD_{n,S}^l \to \cD_{n,T}^l$
is the morphism naturally induced by $f$.
\end{defn}

We will prove that it is actually a DM stack.

\begin{defn}
We define the {\em stack of twisted $n$-jet} of $\cX$ by
\[
 \cJ _n \cX := \coprod _{\chara (k)  \nmid  l} \cJ_n ^l \cX.
\]
\end{defn}

If $\cX$ is of finite type, then $\cJ_n^l \cX$ is empty for sufficiently large $l$ and $\cJ_n \cX$ is in fact
the disjoint sum of only finitely many $\cJ_n^l \cX$.

\begin{lem}
The category $\cJ _n^l \cX$ is a stack.
\end{lem}

\begin{proof}
For an object $\gamma : \cD_{n,S}^l \to \cX$ of $\cJ _n^l \cX$ and for a morphism $f :T \to S$,
we have a``pull-back", $\gamma_T := f' \circ \gamma$ which is unique up to 2-isomorphisms.
Here $f' : \cD_{n,T}^l \to \cD_{n,S}^l $ is the natural morphism induced by $f$. 
Hence $\cJ _n^l \cX$ is a groupoid.

We first show that for two objects 
$\gamma, \gamma' :\cD_{n,S}^l \to \cX$, the functor
\begin{align*}
\Isom (\gamma , \gamma'): (\mathrm{Aff}/S) &\to (\text{Sets}) \\
(T\to S) &\mapsto \Hom _{(\cJ_n^l\cX)(T)}(\gamma_T,\gamma'_T) .
\end{align*}
is a sheaf. Consider a morphism $T \to S$ and an \'etale cover $\coprod T_i \to T$.
Let $T_{ij} := T_i \times _T T_j $. For every object $\alpha$ of $\cD_{n,T}^l$,  
we have the pull-backs $\alpha_i$ and $\alpha_{ij}$ to $\cD_{n,T_i}^l$ and $\cD_{n, T_{ij}}^l  $
respectively. Since $\cX$ is a prestack, the sequence
\begin{align*}
 0 \to \Hom _{\cX(T)}(\gamma_{T}(\alpha), \gamma'_{T}(\alpha)) \to 
 \coprod \Hom_{\cX({T_i})}(\gamma_{T_{i}}(\alpha_{i}), \gamma'_{T_{i}}(\alpha_{i}))\\
  \rightrightarrows \coprod \Hom_{\cX({T_{ij}})}(\gamma_{T_{ij}}(\alpha_{ij}), \gamma'_{T_{ij}}(\alpha_{ij}))
\end{align*}
is exact. Since a morphism of twisted jets is a natural transformation of functors, it
implies that the sequence
\begin{align*}
 0 \to \Hom _{(\cJ_n^l\cX)(T)}(\gamma_{T}, \gamma'_{T}) \to 
\coprod \Hom_{(\cJ_n^l\cX)({T_i})}(\gamma_{T_{i}}, \gamma'_{T_{i}})\\
  \rightrightarrows  \coprod \Hom_{(\cJ_n^l\cX)({T_{ij}})}(\gamma_{T_{ij}}, \gamma'_{T_{ij}})
\end{align*}
is also exact, and the functor $\Isom (\gamma , \gamma')$ is a sheaf.

It remains to show that one can glue objects. Let $\coprod T_i \to T $ be an \'etale cover, let $\gamma_i : \cD _{n,T_i}^l  \to \cX$ be  twisted jets  and 
let $h_{ij}: (\gamma_i)_{T_{ij}} \to (\gamma_j)_{T_{ij}}  $ be a morphism in $(\cJ_n^l\cX)(T_{ij})$.
Assume that they satisfy the cocycle condition. Then for every object $\alpha$ of $\cD_{n,T}^l $,
we can glue the objects $\gamma_i (\alpha_i) $ of $\cX$, because $\cX$ is a stack.  
Therefore we can determine the image of $\alpha$ and obtain a functor $\gamma:\cD _{n,T}^l  \to \cX$
which is clearly representable. Thus we have verified all conditions.
\end{proof}

\subsubsection{}
Let $\cY \to \cX$ be a representable morphism of DM stacks. Then for a twisted jet $\cD _{n,S}^l  \to \cY$,
composing the morphisms, we obtain a twisted jet $\cD _{n,S}^l  \to \cX$.
Thus we have a natural morphism $\cJ_n^l \cY \to \cJ_n^l \cX$.

In \cite{Twistedjet}, we defined a \textit{barely faithful morphism} to be a morphism $f:\cY \to \cX$ 
of DM stacks such that for every object $\xi$ of $\cY$, the map $\Aut (\xi ) \to \Aut (f(\xi))$ is
bijective. Thus all barely faithful morphisms are representable.
 Barely faithful morphisms are stable under base change \cite[Lemma 4.21]{Twistedjet}.

\begin{lem}\label{lem-barelyfaithfuletale}
Let $\cY \to \cX$ be a barely faithful and formally \'{e}tale morphism of DM stacks.
Then the naturally induced diagram
\[
\xymatrix{
\cJ _n ^l \cY \ar[r] \ar[d] & \cJ _n ^l \cX \ar[d] \\
\cY  \ar[r]& \cX
}
\]
is cartesian.
\end{lem}

\begin{proof}
Consider a commutative diagram
\[
\xymatrix{
S \ar[r] \ar[d] & \cY \ar[d] \\
\cD _{n,S} ^l \ar[r] & \cX
}
\]
where the bottom arrow is representable and the left arrow is a natural one.
Then we claim that there exists a unique morphism $\cD _{n,S} ^l  \to \cY$ which fits into
the diagram. The lemma easily follows from it.

Without loss of generality, we can assume that $S$ is connected. Let 
$\cU \subset \cD _{n,S} ^l\times _ {\cX} \cY$ be the connected component containing
the image of $S$. Then the natural morphism $\cU \to \cD _{n,S} ^l $ is barely faithful, formally \'{e}tale and
bijective, hence an isomorphism. It shows our claim. Thus we obtain an equivalence of categories, 
$ \cY \times _{\cX} \cJ^l_n \cX  \cong \cJ_n ^l \cY$.
\end{proof}

For every DM stack $\cX$, there are finite groups $G_i$, schemes $M_i$ with $G_i$-action and 
a morphism $\coprod _i [M_i/G_i] \to \cX$ which is \'{e}tale, surjective and barely faithful.
Hence thanks to Lemma \ref{lem-barelyfaithfuletale},  in proving that $\cJ_n^l \cX$ is a DM stack, we may assume that
$\cX$ is a quotient stack $[M/G]$. 
Let $k' /k$ be the field extension by adding all $l$-th roots of unity for the order $l$ of elements
of $G$ prime to the characteristic of $k$.
Replacing $k$ with $k'$ and $M$ with $M \otimes _k k'$, we
may assume that $\bmu_{l,k}$ is a constant group scheme for 
$l$ such that there is a twisted jet $\cD_{n,S}^l \to [M/G]$.
The action $\bmu_{l} \curvearrowright D_{n,S}$ induces an action $\bmu_l \curvearrowright J_n M$.
On the other hand, 
for each embedding $a:\bmu_l \hookrightarrow G$, $\bmu_l$ acts on $M$ as a subgroup of $G$ and on
$J_n M$. 

\begin{defn}
We define $J^{(a)}_n M$ to be the closed subscheme of $J_n M$ where 
the two actions $\bmu_l \curvearrowright J_n M$ are identical.
\end{defn}

\begin{defn}
We define $\Conj(\bmu_l,G)$ to be a representative set of the conjugacy classes of
embeddings $ \bmu_l \hookrightarrow G $. 
\end{defn}

\begin{prop}\label{lem-verifyDMstack}
For $0 \le n \le \infty$, there is an isomorphism
\[
\cJ _n^l \cX \cong \coprod _{a \in \Conj (\bmu_l,G)} [J_{nl}^{(a)} M/C_a].
\]
Here $C_a$ is the centralizer of $a$. By this isomorphism,
$[J_{nl}^{(a)} M/C_a]$ corresponds to twisted jets $\cD_{n,S}^l \to \cX $
inducing $a : \bmu _l \hookrightarrow G$.
\end{prop}

\begin{proof}
Let $m := nl$. Choose a primitive $l$-th root   $\zeta \in \bmu _l$ of unity. 
 Let $\gamma : \cD_{n,S}^l \to \cX$ be an object  over $S$ of $\cJ _n^l \cX$.
The canonical atlas $D_{m,S} \to \cD_{n,S}^l$ corresponds to the object $\alpha$ 
of  $\cD_{n,S}^l$
\[\xymatrix{
D_{m,S} \times   \bmu_{l} \ar[r]^{\bmu_l\text{-action}} \ar[d]_{\text{trivial $\bmu_l$-torsor}} & D_{m,S} \\
D_{m,S}   & .
}
\]
The morphism 
\[
 \theta:= \zeta \times \zeta^{-1}:  D_{m,S} \times   \bmu_{l}  \to  D_{m,S} \times   \bmu_{l} 
\]
is an automorphism of $\alpha$ \textit{over} $\zeta : D_{m,S} \to D_{m,S}$, whose order is $l$.
Any other object of $\cD _{n,S} ^l$ is a pull-back of $\alpha$ and any automorphism is a pull-back of a power of 
 $\theta$. 
Therefore the twisted jet $\gamma$ is determined by the images of $\alpha$ and $\theta$ in $\cX$. 

Let the  diagram
\[
\xymatrix{
P \ar[r]^h \ar[d]_p & M \\
D_{m,S} &
}
\]
be the object over $D_{m,S}$ of $\cX$ which is the image of $\alpha$ by $\gamma$. Let $\lambda$ be its automorphism
over  $\zeta : D_{m,S}  \to D_{m,S}$ which is the image of $\theta$.
Because $\gamma$ is representable, the order of $\lambda$ is also $l$.
Let $Q:= P \times _{D_{m,S}} S$.
Then $P \to D_{m,S}$ is isomorphic as torsors to $D_{m,k} \times_k Q \to D_{m,k} \times_k S$.
Since we have chosen a primitive $l$-th root $\zeta$, we can identify $\Conj (\bmu_l ,G)$ with
a representative set
 $\Conj^l (G)$ of the conjugacy classes of {\em elements} of order $l$.

\textbf{Claim:} \textit{If $S$ is connected, then there are open and closed subsets 
$Q' \subset Q$ and $P' \subset P$ which, for some $ g   \in \Conj ^l(G)$, are stable under
$C_g$-action and $C_g$-torsors over $S$.}

Take an \'etale cover $T \to S$ such that $Q_T := Q \times _S T$ is isomorphic to the trivial $G$-torsor
 $T \times G \to T$ with a right action. Then the pull-back of the automorphism $\lambda$ is 
a left action of some $g^{-1} \in G$ over each connected component of $T$. 
If necessary, replacing the isomorphism $Q_T \cong T \times G$, we can assume
that the automorphism is given by unique $g^{-1} \in \Conj ^l(G)$. Let $\phi : T \times G \to Q$ be the natural
morphism. Then we see that $\phi (T\times C_g)  \cap \phi(T \times (G \setminus C_g)) = \emptyset$, as follows:
Let $a \in C_g$, $b \in G\setminus C_g$, $x\in T\times C_g$ and $y \in 
T \times (G \setminus C_g)$. If $\phi(x) = \phi(y)$, then 
\[
\phi (x) = \phi(gxg^{-1})=\lambda \phi(x) g^{-1} = \lambda \phi(y) g^{-1}=\phi(gyg^{-1}) \neq \phi(y).
\]
It is a contradiction. Similarly $P$ decomposes also.

Since $(h \circ \lambda) | _ {P'} = h |_ {P'}$ and $ (h \circ g )|_{P'} = (g \circ h)|_{P'} $, 
we have 
\[
h \circ (\zeta \times \id _{Q'} ) =h \circ ( \lambda \circ g  )|_{P'} 
= (g \circ h)|_{P'} .
\]
It means that the morphism $ D_{m,k} \times_k Q' \to M $ corresponds to a morphism $Q' \to (J_{m} M)^{\zeta\circ g^{-1}}$
and that the morphism $Q' \to (J_{m} M)^{\zeta\circ g^{-1}}$ and a $C_g$-torsor $Q' \to S$ determine an object
over $S$ of a quotient stack $[(J_{m} M)^{\zeta\circ g^{-1}}/C_g]$.
Note that $(J_m M)^{\zeta \circ g^{-1}}= J_m^{(a)}M$.
Thus we have a morphism $\cJ_n^l \cX \to \coprod [ J_m^{(a)} M/ C_a ]$.
The inverse morphism can be constructed by following the argument conversely.
\end{proof}

\begin{thm}\label{thm-JnM=DMstack}
Let $\cX$ be a DM stack.
\begin{enumerate}
\item  For  $ n \in \ZZ_{\ge0}$, $\cJ_n^l \cX$ and $\cJ_n \cX$ are DM stacks.
\item If $\cX$ is of finite type (resp.\ smooth), then for $ n \in \ZZ_{\ge0}$,
 then $\cJ _n ^l \cX$ and $\cJ_n \cX$ are also of finite type (resp.\ smooth).
\item For every $m \ge n$, the natural projection $\cJ_m \cX \to \cJ_n \cX$ is an affine morphism.
\end{enumerate}
\end{thm}

\begin{proof}
1: There is an \'{e}tale, surjective and barely faithful morphism 
$\coprod _i [M/G_i] \to \cX$ such that each $M_i$ is a scheme and $G_i$ is a finite group.
From  Lemma \ref{lem-barelyfaithfuletale}, Proposition \ref{lem-verifyDMstack} and \cite[Lemme 4.3.3]{LMB},
the morphism $\cJ ^l _n \cX \to \cX$ is representable. From \cite[Proposition 4.5]{LMB}, 
$\cJ^l_n \cX$ is a DM stack.  $\cJ _n \cX$ is also a DM stack.
The morphism  $\cJ ^l _n \cX \to \cX$ is also separated and so is $\cJ ^l _n \cX$.

2 and 3: These also result from Lemma \ref{lem-barelyfaithfuletale} and Proposition \ref{lem-verifyDMstack}.
\end{proof}

In general, a  projective system 
$\{ \cX_i,\rho_i : \cX_{i+1} \to \cX _i \}_{i \ge 0}$ of DM stacks such
that every $\rho _i$ is representable and affine, 
there exists a projective limit
$\displaystyle \cX _{\infty} = \lim _{\longleftarrow} \cX _i$. 
In fact, 
for each $i$, there is an $\cO_{\cX_0}$-algebra $\cA_i$ such that
$\cX_i \cong \Spec  \cA_i$ (see \cite[\S 14.2]{LMB}) and the  $\cA_i$'s constitute an inductive system.
We can see that
  $ \displaystyle \cX_\infty := \Spec ( \lim _{\longrightarrow} \cA_i)$ is
  the projective limit of the given projective system.

From Theorem \ref{thm-JnM=DMstack},
 the projective system $\{\cJ_n \cX\}_n$ (resp.\ $\{ \cJ^l_n \cX \}_n$) has the projective limit 
 \[
  \cJ _\infty \cX:= \lim _ {\longleftarrow} \cJ_n \cX \ 
  ({resp.\ }\  \cJ^l _\infty \cX:= \lim _ {\longleftarrow} \cJ^l_n \cX).
 \]
Then the point set $| \cJ _\infty \cX |$ is
identified with the set of the equivalence classes of the twisted arcs $\cD_{\infty ,K}^l \to \cX$ with respect to
the following equivalent relation: Let $\gamma _i : \cD_{\infty , K_i}^l \to \cX$, $i=1,2$, be twisted
arcs. If for a field $K_3 \supset K_1,K_2$ and natural morphisms $\cD _{\infty ,K_3}^l \to 
\cD _{\infty ,K_1}^l, \cD _{\infty ,K_2}^l$, the diagram
\[
\xymatrix{
\cD_{\infty,K_3}^l \ar[r] \ar[d] & \cD_{\infty , K_1}^l \ar[d] ^{\gamma_1} \\
\cD_{\infty , K_2}^l \ar[r]_{\gamma_2} & \cX 
}
\]
is commutative, then $\gamma _1$ and $\gamma_2$ are equivalent.

\begin{rem}
For two stacks $\cX$ and $\cY$, we can define a Hom-stack $\mathcal{H}om(\cX,\cY)$ which parameterizes
morphisms from $\cX$ to $\cY$, and its substack $\mathcal{H}om^{rep}(\cX,\cY)$ which parameterizes
representable morphisms. 
Olsson \cite{Olsson-hom-stack} proved that 
if $\cX$ and $\cY$ are Deligne-Mumford stacks satisfying certain conditions, then
$\mathcal{H}om(\cX,\cY)$ is a Deligne-Mumford stack.
 and $\mathcal{H}om^{rep}(\cX,\cY)$ is its open substack. 
Then, Aoki \cite{aoki-hom-stack} proved that $\mathcal{H}om(\cX,\cY)$ is an Artin stack if
$\cX$ and $\cY$ are Artin stacks satisfying certain conditions. 
The stack $ \cJ_n^l \cX $ of twisted $n$-jets of order $l$ ($n< \infty$) is
identical with $\mathcal{H}om^{rep} (\cD_n^l,\cX)$.
\end{rem}

\subsubsection{Inertia stack}

\begin{defn}
To each DM stack $\cX$, we associate the \textit{inertia stack} $I\cX$ defined as follows;
an object of $I\cX$ is a pair $(x , \alpha)$ with $x$ an object of $\cX$  and $\alpha \in \Aut (x)$ and 
a morphism  $(x,\alpha) \to (y,\beta)$  in $I\cX$ is a morphism $\phi :x \to y$ in $\cX$ 
with $\phi \alpha = \beta \phi$. 
\end{defn}

It is known that $I\cX$ is isomorphic to $\cX \times _{\Delta , \cX \times \cX , \Delta} \cX$, where
$\Delta: \cX \to \cX \times \cX$ is the diagonal morphism. Then the forgetting morphism 
$I\cX \to \cX$ is isomorphic to the first projection
$\cX \times _{\Delta , \cX \times \cX , \Delta} \cX \to \cX$.
Since we have supposed that $\cX$ is separated, 
the diagonal morphism is finite and unramified. 
Hence 
the forgetting morphism $I\cX \to \cX$ is so as well.

\begin{defn}
Let $l$ be a positive integer prime to $\chara (k)$.
We define $I^l \cX \subset I\cX$ to be the open and closed substack of objects $(x,\alpha)$ such that
the order of $\alpha$ is $l$.
\end{defn}

\begin{prop}
Assume that $k$ contains all $l$-th roots of unity.
Then for each choice of a  primitive $l$-th root $\zeta$ of unity,
there is a natural isomorphism $\cJ^l _0 \cX \cong I^l \cX$.
\end{prop}

\begin{proof}
The assertion follows from the fact that
 giving a representable morphism $\cD _0 ^l \times S \to \cX$
is equivalent to giving an object $x$ over $S$ of $\cX$ and 
an embedding $\bmu_l \hookrightarrow \Aut (x)$, which is equivalent to giving the image of $\zeta \in \bmu_l$. 
\end{proof}

The inertia stack is the algebraic counterpart of the twisted sector of an analytic orbifold, which
was used to define the orbifold cohomology in \cite{CR}. 
Since for $k=\CC$, there is a canonical choice $\exp(2\pi \sqrt{-1} /l)$ of a primitive $l$-th root of unity,
we have a natural isomorphism $\cJ_0 \cX \cong I \cX$.

\subsection{Morphism of stacks of twisted jets}

As we saw above,
for a representable morphism $\cY \to \cX$ of DM stacks, we have a naturally induced morphism
$\cJ _ n ^l \cY \to \cJ _ n ^l \cX$. 
For a morphism $\cX \to X$ from a DM stack to an algebraic variety, we can associate 
a morphism $\cJ _ n ^l \cX \to J_n X$ as follows: For a twisted jet $\cD_{n,S}^l \to \cX$, consider
the composition $\cD_{n,S}^l \to \cX \to X$. From the universality of the coarse moduli space, it 
uniquely factors as
$\cD_{n,S}^l  \to D_{n,S} \to X$ up to isomorphism of $D_{n,S}$. 
Putting the condition that $\cD_{n,S}^l=[D_{nl,S}/\bmu_{l,k}]  \to D_{n,S}$ is defined by
$t \mapsto t^l$,
we obtain a unique jet on $X$.

We generalize these to a general (not necessarily representable) morphism of DM stacks.

\begin{prop}\label{lem-morphismoftwistedjetstacks}
Let $f:\cY \to \cX$ be a morphism of DM stacks.
Then for $0 \le n \le \infty$, we have a natural morphism $f_n:\cJ _ n  \cY \to \cJ _ n \cX$.
\end{prop}

\begin{proof}
We may assume $n < \infty$.
Let $\cD_{n,S}^l  \to \cY$ be an object over a scheme $S$ of $\cJ _n ^l \cY$.
Then the composite $\cD_{n,S}^l  \to \cY \to \cX$ is not in general representable.
Take the canonical decomposition $\cD_{n,S} ^l  \to \cE \to \cX$ as in the following lemma.
If $S$ is connected, $\cE$ must be isomorphic to $\cD_{n,S} ^{l'} $ for some divisor $l'$
of  $ l$.
We choose the isomorphism $\cE \cong \cD_{n,S}^{l'}$ so that 
the morphism $  \cD_{n,S} ^l \to \cE \cong \cD_n^{l'}$ induces a morphism $D_{nl,S}  \to D_{nl',S}$ of atlases defined by $t \mapsto t^{l/l'}$.
Thus we obtain a twisted jet $ \cD_{n,S} ^{l'}  \to \cX $ over $\cX$ and 
a morphism $\cJ _ n  \cY \to \cJ _ n \cX$.
\end{proof}

\begin{lem}[Canonical factorization]
Let $f:\cW \to \cV$ be a morphism of DM stacks. Then there are a DM stack $\cW'$ and morphisms
$g:\cW \to \cW'$ and $h:\cW' \to \cV$ such that
\begin{enumerate}
\item $f = h \circ g$,  
\item $h$ is representable, and
\item (universality) 
Let $g':\cW \to \cZ$ and $h':\cZ \to \cV$ be morphisms of DM stacks.
If $h'$ is representable and if the diagram of solid arrows
\[\xymatrix{
\cW \ar[r]^{g'} \ar[d]_{g} & \cZ \ar[d]^{h'} \\
\cW'\ar[r]_{h} \ar@{-->}[ur]^{i} & \cV 
}
\]
is commutative, then there is a morphism $i:\cW' \to \cZ$ making the whole diagram commutative, which is 
unique up to unique 2-isomorphism.
\end{enumerate}
\end{lem}

\begin{proof}
Let $M \to \cV$ be an atlas. It induces a groupoid space $M \times _{\cV} M \rightrightarrows  M$.
The stack associated to the groupoid space is canonically identified with $\cV$.
Let $\cW _M := \cW \times _{\cV} M$.
Similarly we can obtain some structure $\cW_M \times _{\cW} \cW_M \rightrightarrows  \cW_m$ like
a groupoid space, but each object is not a variety but a DM stack.
Taking the coarse moduli space of the objects, we obtain a groupoid space 
$\overline{\cW_M \times _{\cW} \cW_M} \rightrightarrows  \overline{\cW_M}$.
Note that the existence of the coarse moduli space was proved by Keel and Mori \cite{Keel-Mori}.
Let $N \to \cW_M$ be an atlas. Then the composite $N \to \cW_M \to \cW$ is also an atlas.
So we have another groupoid space $N \times _{\cW} N \rightrightarrows N$.
Now there are naturally induced 
morphisms of the groupoid spaces
\[
(N \times _{\cW} N \rightrightarrows N)\to 
(\overline{\cW_M \times _{\cW} \cW_M} \rightrightarrows  \overline{\cW_M})
 \to(M \times _{\cV} M \rightrightarrows  M),
\]
and the corresponding morphisms of DM stacks
\[
\cW \to \cW'  \to \cV
\]
where $\cW'$ is the stack associated to $\overline{\cW_M \times _{\cW} \cW_M} \rightrightarrows  \overline{\cW_m}$.
They clearly satisfy Conditions 1 and 2. 

Suppose that there are morphisms $g' : \cW \to \cZ$ and $h' : \cZ \to \cV$ as in Condition 3.
Taking the fiber products $\times _\cV M$, we obtain the following diagram of solid arrows
\[\xymatrix{
\cW_M \ar[r] \ar[d] & \cZ _M \ar[d] \\
\overline{\cW_M} \ar[r] \ar@{-->}[ur]^{i_M} & M .
}
\]
Since $h'$ is representable, $\cZ_M$ is isomorphic to an algebraic space. 
Now from the universality of the coarse moduli space, there is a morphism 
$i_M $  making the whole diagram commutative, which is unique up to unique 2-isomorphism. 
It implies the last condition.
\end{proof}

Because of the universality, the morphisms $g$ and $h$ are uniquely determined.
We call it the \textit{canonical factorization} of $f$.
For each point $y \in \cY$ and its image $x\in \cX$, we have a homomorphism of the
automorphism groups $\phi : \Aut (y) \to \Aut (x)$. Then the canonical factorization corresponds to 
the factorization of $\phi$, 
\[
\Aut (y) \twoheadrightarrow  \Im (\phi) \hookrightarrow  \Aut (x).
\]
The canonical factorization is a generalization of the coarse moduli space.
Indeed, in the lemma, if $\cV$ is an algebraic space, then $\cW'$ is the coarse moduli space of
$\cW$.

%%%%%%%%%%%%%%%%%%%%%%%%%%%%%%%%%%%%%%%%%%%%%%%%%%%%%%%%%%%%%%%%%%%%%%%%%%%%%%%%%%%%%%%%%%%%
\section{Motivic integration}

\subsection{Convergent stacks}

In this subsection, we construct the semiring in which integrals take values.

\begin{defn}\label{def-convergent-stack}
A {\em convergent stack} is the pair $(\cX ,\alpha)$ 
of a DM stack $\cX$ and a function 
\[
\alpha : \{\text{connected component of $\cX$}\} \to \ZZ
\]
 such that 
\begin{enumerate}
\item there are at most countably many connected components,
\item all connected components are of finite type, and
\item for every $ m \in \ZZ$, there are at most finitely many connected components $\cV$
with $\dim \cV  + \alpha (\cV) > m$.
\end{enumerate}
We say that $\cX$ is the {\em underlying stack} of $(\cX,\alpha)$. 
\end{defn}

We abbreviate $(\cX , \alpha)$ as $\cX$, if it  causes no confusion.
A DM stack of finite type $\cX$ is identified with the convergent stack
$(\cX,0)$.

A {\em  morphism} $f: (\cY , \beta ) \to (\cX ,\alpha) $ of convergent stacks is a morphism 
$f:\cY \to \cX$ of
the underlying stacks with $ \beta = \alpha \circ f$.
A morphism of convergent stacks is called an {\em isomorphism} if
it is an isomorphism of the underlying stacks.
We say that convergent stacks $\cX $ and $\cY$ are {\em  isomorphic} (write $\cX \cong \cY$)
if there is an isomorphism between them.

If $(\cX,\alpha)$ is a convergent stack and $\cX' $ is a locally closed substack of $\cX$,
and if $\iota : \cX' \hookrightarrow \cX$ is the inclusion map, 
then $(\cX', \alpha \circ \iota)$ is a convergent stack 
We call this a {\em convergent substack} of $\cX$.

In the category of convergent stacks, we have the disjoint union and the product of two objects:
Let $(\cX ,\alpha)$ and $(\cY , \beta )$ be convergent stacks.
The disjoint union $(\cX ,\alpha)\sqcup (\cY , \beta )$
 of them is the convergent stack $( \cX \sqcup  \cY , \gamma )$ such that
 $\gamma |_\cX= \alpha$ and $\gamma |_\cY = \beta$.
 The product  $(\cX ,\alpha)\times (\cY , \beta )$ is 
 the convergent stack $( \cX \times  \cY , \gamma )$ such that
 for a connected component $\cV \subset \cX \times \cY$, 
 $\gamma (\cV) = \alpha ( p_1 (\cV) ) + \beta (p_2 (\cV))$.

\begin{defn}
For a convergent stack $\cX=(\cX , \alpha)$,
 we define the {\em dimension} of $\cX$, denoted $\dim \cX$, to be 
\[
\max \{\dim \cV + \alpha (\cV) | \cV \subset \cX \text{ connected component} \}.
\]
By convention, we put $\dim \emptyset = - \infty $.
\end{defn}

\begin{defn}\label{def-equivalence-relation}
Let $\fR '$ be the set of the isomorphism classes of convergent stacks.
For each $n \in \ZZ$, we define $\sim _n$ to be the strongest equivalence relation of $\fR'$
satisfying the following basic relations:
\begin{enumerate}
\item If $\cX$ and $\cY$ are convergent stacks  with $\dim \cY < n$, then $\cX \sim _n \cX \sqcup \cY$.
\item If $\cX \subset \cY$ is a convergent closed substack, then
$ \cY \sim _n (\cY \setminus \cX) \sqcup \cX $. (Here $\sqcup$ is not the disjoint union in $\cY$.)
\item Let $(\cY,\beta)$ and $(\cX , \alpha)$ be convergent stacks and
$f : \cY \to \cX$ a representable morphism of underlying stacks.
If every geometric point $y$ of $\cY$, $f^{-1}f(y)$ is isomorphic to an affine space of
dimension $ \alpha(f(y))- \beta(y) $, then $(\cY , \beta) \sim_n (\cX, \alpha)$.  
\end{enumerate}
Namely $\cX \sim_n \cY$ if and only if $\cX$ and $\cY$ can be connected by
finitely many basic relations above.
We define an equivalence relation $ \sim $ of $\fR'$ as follows:
For $a,b \in \fR'$, $a \sim b$ if and only if $a \sim _n b$ for all $n \in \ZZ$.
\end{defn}

For example, we have
\begin{align*}
 (\Spec k, 0) &\sim (\AA^1_k ,-1) \\
& \sim (\AA^1_k \setminus \{ 0 \} ,-1) \sqcup (\AA^1 , -2) \\
& \sim  \coprod _{i \in \NN} (\AA^1_k \setminus \{0\}, -i) .
\end{align*}

\begin{defn}
We define $\fR $ to be $\fR'$ modulo $\sim$.
For a convergent stack $\cX$, we denote by $\{\cX\}$ the equivalence class of $\cX$. 
\end{defn}

\begin{lem}
Let $\cX$ and $\cY$ be convergent stacks.
If $\cX \sim \cY$, then $\dim \cX = \dim \cY$.
\end{lem}

\begin{proof}
By definition, for all $n $, $ \cX $ and $\cY$ are connected by finitely many
basic relations in Definition \ref{def-equivalence-relation}. 
For $n \ll 0$, these relations preserve the dimension. 
Hence $\dim \cX = \dim \cY$.
\end{proof}

From the lemma, we have a map
\[
 \dim : \fR \to \ZZ \cup \{- \infty \},\ \{\cX\} \mapsto \dim \cX .
\]

\begin{lem}
Let $\cX_1$, $\cX_2$, $\cY_1$ and $\cY_2$ be convergent stacks.
If $\cX_1 \sim \cX_2$ and $\cY_1 \sim \cY_2$, 
then $ \cX _1 \sqcup \cY _1 \sim \cX_2 \sqcup \cY_2$ and 
$ \cX _1 \times \cY _1 \sim \cX_2 \times \cY_2$.
\end{lem}

\begin{proof}
The assertion $ \cX _1 \sqcup \cY _1 \sim \cX_2 \sqcup \cY_2$ is clear.

If $\cX_1 \sim _n \cX_2$, then 
$\cX_1$ and $\cX_2$ are connected by the basic relations in Definition \ref{def-equivalence-relation}.
Therefore $\cX_1 \times \cY \sim _ {n + \dim \cY} \cX_2 \times \cY$ for any convergent stack $\cY$.
This implies that if $\cX_1 \sim \cX_2$, then $ \cX_1 \times \cY \sim  \cX_2 \times \cY$.
This proves the second assertion.
\end{proof}

We define a commutative semiring structure on $\fR$ by $\{ \cX \} + \{\cY\} := \{\cX \sqcup \cY\}$ and
$\{\cX\}\{\cY\}:= \{\cX \times \cY\}$.
The element $\{ \emptyset\}$ is the unit of addition and the element $\{ \Spec k \} = \{(\Spec k,0)\}$
is the unit of product. 
We denote the element $\{\AA^1_k\} = \{ (\Spec k,1)\}$ by $\LL$.
This is clearly invertible ($\LL^{-1} = \{ (\Spec k, -1) \}$).

Let $I$ be a countable set and let $\cX_i$, $i \in I$, be convergent stacks.
If for every $m \in \ZZ$, there are at most finitely many $i \in I$ with $\dim \cX_i > m$,
then  the disjoint union $ \coprod _{i \in I} \cX_i $ is also a convergent stack.
In other words, if $x_i \in \fR$, $i \in I$ and if there are finitely many 
$i \in I$ with $\dim x_i > m$, then  the infinite sum $\sum _{i \in I} x_i $
is  defined.
Note that  $\sum _{i \in I}x_i$ is  independent of the
order of $I$. 
Moreover the following holds:
If $I_j \subset I$ are subsets with $I= \coprod_j I_j$,
then we have
\[
 \sum_{i \in I} x_i =  \sum _{j} \sum _{i \in I_j} x_i.
\]

\begin{defn}
We define
$\bar \fR := \fR \cup \{\infty\}$.
\end{defn}

 This is still a semigroup by 
$x + \infty = \infty$ for any $x \in \bar \fR $.
Let $I$ be a countable set and $x_i \in \bar \fR$, $i \in I$.
Either if there exists $i \in I$ with $x_i = \infty$ or if
 for some $m \in \ZZ$, there exist infinitely many $i$ with $\dim x_i > m$,
then we define $\sum_{i \in I} x_i := \infty $.
Thus in the semigroup $\bar \fR$, the sum of arbitrary countable collection of elements is defined.

We can see the following basic properties of the map $\dim$.

\begin{lem}\label{lemma-properties-dim}
\begin{enumerate}
\item For $x,y\in \fR $, $\dim (x \times y) = \dim x + \dim y$.
\item Let $I$ be a countable set and $x_i \in \fR$, $i \in I$. 
If $ \sum _{i\in I}x_i \neq \infty$, then 
 $\dim \sum _{i\in I}x_i = \max \{\dim x_i | i \in I\}$.
\end{enumerate}
\end{lem}

\begin{rem}
Originally, the motivic integration was defined in a ring deriving from the Grothendieck ring of
varieties, in precise, a completion of a localization of the Grothendieck ring (see \cite{germs}).
In the Grothendieck ring, there are  negative elements, that is, there is the inverse of 
any class of varieties. This causes some problems.
For example, the Grothendieck ring of varieties is not a domain \cite{Poonen-not-domain}. 
However we do not really need negative elements to construct the integration theory. 
Thus it is natural to replace the Grothendieck ring with the Grothendieck {\em semiring}.

However we do not superficially use the Grothendieck semiring. 
Instead, we consider a sufficiently large set (the set of the isomorphism classes of convergent stacks or
 spaces) and take its quotient by relations needed below.
Advantages of this approach are that each element of the obtained semiring
 corresponds to a ``geometric'' object
and that we do not have to use algebraic operations on semirings like the localization and
the completion.
\end{rem}

\subsection{Convergent spaces and coarse moduli spaces}

A {\em convergent space} is defined to be a convergent stack $(\cX, \alpha )$ such that
the underlying stack $\cX$  is an algebraic space.

\begin{defn}\label{def-equivalence-relation2}
Let $\fS '$ be the set of the isomorphism classes of convergent spaces.
For each $n \in \ZZ$, we define $\simeq  _n$ to be the strongest equivalence relation of $\fS'$
satisfying the following basic relations:
\begin{enumerate}
\item If $f : Y \to X$ is a morphism of convergent spaces which is 
finite, surjective and universally injective, then $X \simeq _n Y$.
\item If $X$ and $Y$ are convergent spaces  with $\dim Y < n$, then $X \simeq _n X \sqcup Y$.
\item If $X \subset Y$ is a convergent closed subspace, then
$ Y \simeq _n (Y \setminus X) \sqcup X $. 
\item Let $(Y,\beta)$ and $(X , \alpha)$ be convergent spaces and
$f : Y \to X$ a morphism of underlying algebraic spaces.
If every geometric point $y \in Y (K)$, 
there is a finite, surjective and universally injective morphism 
\[
 \AA^{\alpha(f(y))- \beta(y) }_K /G \to f^{-1}f(y)
\]
for some finite group action $G \curvearrowright \AA^{\alpha(f(y))- \beta(y) }_K $,
then $(Y , \beta) \simeq _n (X, \alpha)$.  
\end{enumerate}
We define an equivalence relation $ \simeq $ of $\fR'$ as follows:
For $a,b \in \fR'$, $a \simeq b$ if and only if $a \simeq _n b$ for all $n \in \ZZ$.
\end{defn}

A finite, surjective and universally injective morphism
induces an equivalence of \'etale sites \cite[IX. Th\'eor\`eme 4.10]{SGA1} and
an isomorphism of \'etale cohomology.
We will use this fact in the following subsection.  

\begin{defn}
We define $\fS$ to be $\fS'$ modulo $\simeq$. 
\end{defn}

The set $\fS$ has a semiring structure as $\fR$ does.

For a DM stack $\cX$, we denote by $\bar \cX$ its coarse moduli space.
For a convergent stack $\cX = (\cX , \alpha)$,
we can give a natural convergent space structure to the coarse moduli space $\bar \cX$.
We denote this convergent space also by $\bar \cX$.    

\begin{prop}
There is a semiring homomorphism 
\[
\fR \to \fS,\ \{\cX\} \mapsto \{\bar \cX\}.
\]
\end{prop}

\begin{proof}
First, to show that there is a map of sets, 
we have to show that for each $n$, if $\cX$ and $\cY$ are
in one of the basic relations in Definition \ref{def-equivalence-relation},
then $\bar \cX \simeq _n \bar \cY $. 

Concerning the first basic relation, this is clear.
We next consider the second one.
If $\cX \subset \cY$ is a convergent closed substack and $X \subset \bar \cY$
is a convergent closed subspace which is the image of $\cY$, then 
 the natural morphism $\bar \cX \to X$ is finite, surjective and universally injective.
Thus $\bar \cY \simeq _n \overline{\cY \setminus \cX} \sqcup \bar \cY$.

Let $f:\cY \to \cX$ be a morphism of DM stacks and 
$g : Y \to X$ the corresponding morphism of the coarse moduli spaces.
Let  $x :\Spec K \to \cX$ be a geometric point and $\bar x : \Spec K \to X$ the corresponding point of the coarse moduli space. Suppose that
$f^{-1}(x) \cong \AA_K^d$.
Then there is a natural morphism $\AA^d_K/\Aut (x) \to g^{-1}(\bar x)$
which is finite, surjective and universally injective. 
This implies that if $\cX$ and $\cY$ is in the last basic 
relation in Definition \ref{def-equivalence-relation},
then $\bar \cX$ and $\bar \cY$ are in the last basic 
relation in Definition \ref{def-equivalence-relation2}.

Thus we have a map $\fR \to \fS,\ \{\cX\} \mapsto \{\bar \cX\}$.

Next, to show that this map is a semiring homomorphism,
we have to show that 
\[
\{\bar \cX\} \{\bar \cY\}:= \{\bar \cX \times \bar \cY\} = \{ \overline{\cX \times \cY} \}.
\]
There is a natural morphism $\overline{\cX \times \cY} \to \bar \cX \times \bar \cY$.
This is finite, surjective and universally injective.
It follows that $\{\bar \cX \times \bar \cY\} = \{ \overline{\cX \times \cY} \}$.
\end{proof}

\subsection{Cohomology realization}\label{subsec-cohomology-realization}

To control huge semirings $\fR$ and $\fS$, it is useful to consider homomorphisms
from these semirings to smaller and easier ones. The cohomology theory helps to construct such homomorphisms. 

\subsubsection{The Grothendieck ring of an abelian category}

Let $\cC$ be an abelian category.
Its {\em Grothendieck ring} $K_0(\cC)$ is
defined to be the free abelian group generated by isomorphism classes $[M]$
of objects $M$ of $\cC$
modulo the following relation: If there is a short exact sequence
\[
 0 \to M_1 \to M_2 \to M_3 \to 0,
\]
then $[M_2] = [M_1]+[M_3]$. 
It is easy to see that given a sequence $0 = M_0 \subset M_1 \subset \dots \subset M_n =M$, then we have an equation 
\[
 [M] = \sum_{i =0}^{n-1} [M_{i+1}/M_i] .
\]
Every object $M$ of $\cC$ has a Jordan-H\"older sequence
$0 = M_0 \subset M_1 \subset \dots \subset M_n =M$, that is, 
all $M_{i+1}/M_i$ are simple objects. Then the {\em semisimplification} $M^{ss}$ is by definition
the associated graded $\bigoplus_i M_{i+1}/M_i$ of the sequence. 
It is a semisimple object and its isomorphism class depends only on $M$. 
We have $[M]=[M^{ss}]$.
If $M_1$ and $M_2$ are semisimple and $[M_1]=[M_2] \in K_0 (\cC)$, then 
$M_1$ and $M_2$ are isomorphic. If $M_1$ and $M_2$ are arbitrary objects
 and $[M_1]=[M_2] \in K_0 (\cC)$, then 
$M_1^{ss} \cong M_2^{ss}$.

\subsubsection{Mixed Hodge structures}

Now we consider the case where  $\cC$ is the category of (rational) mixed Hodge structures, denoted $MHS$.
Since every mixed Hodge structure has, by definition, a weight filtration
whose associated graded is a pure Hodge structure, 
the Grothendieck ring $K_0 (MHS)$ is, in fact, generated by the classes of pure Hodge structures. For a variety $X$ over $\CC$, we define 
\[
 \chi_h(X):= \sum_i (-1)^{i} [H^i_c(X ,\QQ)] \in K_0(MHS).
\]
For a variety $X$ and its closed subvariety $V$, since
there is the localization sequence
\[
\cdots \to H_c^i(X \setminus V) \to H_c^i(X) \to H_c^i(V) \to H_c^{i+1}(X \setminus V) \to \cdots ,
\]
 we have $\chi_h (X)= \chi_h (X \setminus V) + \chi _h(V)$. 
We denote  $[\QQ(-1) ] = \chi_h (\AA^1_k)$ by $\LL$.

\subsubsection{$p$-adic Galois representations}

Suppose that $k$ is a finite field.
Let $p$ be a prime number different from the characteristic of $k$ and
 $\QQ_p$ be the $p$-adic field.
Then the compact-supported 
$p$-adic cohomology groups $H^i_c(X \otimes \bar k, \QQ_p)$ are finite dimensional $\QQ_p$-vector spaces with continuous actions
of the absolute Galois group $\G_k = \Gal(\bar k /k)$.

Let $V$ be an arbitrary Galois representation, that is, a
 finite dimensional $\QQ_p$-vector space with continuous $\G_k$-action. 
We say that $V$ is {\em pure of weight $i \in \QQ$}  if 
every eigenvalue $\alpha \in \bar \QQ_p$ of the Frobenius action on $V$ is algebraic and 
all complex conjugates of $\alpha$ have absolute value $p^{- i/2}$. We say that 
$V$ is {\em mixed (of weight $\le i_n, \ge i_0$)}  if there is a filtration 
$0 = W_{i_0} \subset W_{i_1} \subset \cdots \subset W_{i_n} = V $, $(i_j \in \QQ)$ such that
$W_{i_{j}}/W_{i_{j-1}}$ is pure of weight $i_j$. 
(We admit rational weights for later use.)

Let $MR(\G_k,\QQ_p)$ be the category of mixed Galois representations. 
This category is abelian. For a variety $X$, $ H_c^{i}(X \otimes _k \bar k, \QQ_p) $
is mixed \cite{Del-Weil-II}. Define
\[
\chi_p(X):= \sum _i (-1)^{i}[H_c^{i}(X \otimes _k \bar k, \QQ_p)] \in K_0(MR(\G_k, \QQ_p)).
\]
Then this invariant have the same properties as $\chi_h$, namely, 
$\chi_p(X)=\chi_p (X \setminus V) + \chi_p(V)$ for a closed subvariety $V \subset X$,
and $\chi_p (X \times Y) = \chi_p (X) \chi _p (Y)$. 

\subsubsection{Completions of Grothendieck rings}

Let $F_m K_0 (MHS) \subset K_0 (MHS)$ be the subgroup generated by the elements
$[H]$ with $H$ of weight $\le - 2m$. 
 We define the completion 
\[
 \hat K_0 (MHS) := \lim _{\longleftarrow} K_0 (MHS)/F_m K_0 (MHS).
\]
Since $F_m K_0 (MHS)\cdot F_m K_0 (MHS) \subset F_{m+n}K_0 (MHS)$,
$\hat K_0 (MHS)$ has a ring structure.

\begin{lem}
The natural map $K_0 (MHS) \to \hat K_0(MHS)$ is injective.
\end{lem}

\begin{proof}
Let $\alpha = \sum_i n_i [H_i] \in K_0 (MHS)$.
Then if we put $\alpha _w := \sum n_i[\Gr _W^w H_i]$, we have that $\alpha = \sum _w \alpha _w$
and that $\alpha =0$ if and only if $\alpha _w = 0$ for every $w$. 
Moreover $\alpha \in F_m$ if and only if $\alpha _w =0$ for every $w > -2m$.
Now we can see that $\bigcap_m F_m K_0 (MHS) = \{0\}$ and the completion map is injective.
\end{proof} 

Next we define a completion of $K_0(MR(\G_k,\QQ_p))$.
For each integer $m$, we define $F_m K_0(MR(\G_k,\QQ_p)) \subset K_0(MR(\G_k, \QQ_p))$ to be the subgroup
generated by the elements $[V]$ with $V$ of weight $\le -2m$.
We define 
\[
 \hat K_0(MR(\G_k,\QQ_p)) := \lim _{\longleftarrow} K_0(MR(\G_k,\QQ_p)) / F_m K_0(MR(\G_k,\QQ_p)) .
\]
This is also a ring.
The natural map $ K_0(MR(\G_k,\QQ_p)) \to \hat K_0(MR(\G_k,\QQ_p)) $ 
is also injective.

\subsubsection{Maps from $\fR$ and $\fS$}

When $k= \CC$, for a convergent space $ X= (X,\alpha)$, 
we define 
\[
 \chi _h (X ):= \sum_{V \subset X} \chi _h (V) \LL ^{\alpha (V)} \in \hat K_0 (MHS).
\]
Here the sum runs over all connected components $V$ of $X$.
By definition, for every $m \in \ZZ$, there are at most finitely many
$V$ with $\dim V + \alpha (V) >m$.
Since $\chi _h (V)$ is of weight $\le 2 \dim V$, 
for every $m \in \ZZ$, there are at most finitely many
$V$ with 
$\chi _h (V) \LL ^{\alpha (V)} \notin F_ {- m} K_0 (MHS)$.
Therefore $\chi_h (X)$ is well-defined.

When $k$ is a finite field, for a convergent space $X$,
we similarly define 
\[
\chi_p (X) := \sum_{V \subset X} \chi _p (V) \LL ^{\alpha (V)}  \in \hat K_0(MR(\G_k,\QQ_p)).
\]

\begin{prop}
When $k=\CC$, we have a semiring homomorphism 
\[
\fS \to \hat  K_0 (MHS),\ \{X\} \mapsto \chi_h (X).
\]
When $k$ is a finite field, we have a semiring homomorphism
\[
\fS \to \hat  K_0 (MR(\G_k,\QQ_p)),\ \{X\} \mapsto \chi_p (X).
\]
\end{prop}

\begin{proof}
If the maps are well-defined, then 
these are obviously semiring homomorphisms.
Now it suffices to show that for every $n \in \ZZ$,
the basic relations in Definition \ref{def-equivalence-relation2}
preserves $\chi_h (X)$ modulo $F_{-n} \hat K_0 (MHS)$
or $\chi_p(X)$ modulo $F_{-n} \hat K_0 (MR(\G_k,\QQ_p))$.
 Here $F_{-n} \hat K_0 (MHS)$  and $F_{-n} \hat K_0 (MR(\G_k,\QQ_p))$ are the completions 
 of $F_{-n} K_0 (MHS)$ and $F_{-n}  K_0 (MR(\G_k,\QQ_p))$.
The second basic relation clearly does.

The basic relations except for the second one, in fact, preserves
$\chi_h (X)$ and $\chi_p (X)$ even in $\hat K_0 (MHS)$ and 
 $\hat K_0 (MR(\G_k,\QQ_p))$. Only the fourth one is not trivial.
This results from Lemma \ref{lemma-characteristic}
and the fact that a finite, surjective and universally injective morphism is an 
homeomorphism both in analytic topology and in \'etale topology 
(for \'etale topology, see \cite[IX. Th\'eor\`eme 4.10]{SGA1}).
\end{proof}

\begin{lem}\label{lemma-characteristic}
Suppose that $k$ is either $\CC$ or a finite field.
 Let $f: Y \to X$ be a morphism of varieties and $d$ a positive integer.
Suppose that for every geometric point $x\in X(K)$, the fiber $f^{-1}(x)$ 
has the same  compact-supported cohomology as $\AA_K^d$.
Then we have
\begin{align*}
 \chi _h( Y )&= \chi_h (X) \LL^d \ (k= \CC) \\
  \chi _p( Y )&= \chi_p (X) \LL^d \ (k \ 
 \text{finite field}).
\end{align*}
\end{lem}

\begin{proof}
We discuss $\chi_h$ and $\chi_p$ together and write $\QQ_p$ as $\QQ$.
There is a spectral sequence
\[
E_2^{i,j} = H_c^i(X, R^jf_! \QQ) \Rightarrow H_c^{i+j}(Y,\QQ)
\]
From Grothendieck's generic flatness, we may assume that $f$ is flat.
Then from the assumption, we have
\[
R^jf_! \QQ = 
\begin{cases}
\QQ (-d) &  (j=2d) \\
0 & (j \ne 2d) ,
\end{cases}
\]
and the spectral sequence degenerates. Hence we have an isomorphism
\[
 H_c^{i+2d} (Y,\QQ) \cong H_c^i (X,\QQ) \otimes \QQ(-d),
\]
which implies the lemma.
\end{proof}

Composing maps, we obtain  semiring homomorphisms
\begin{align*}
 \fR \to \hat K_0 (MHS),\  \cX \mapsto \chi _h (\bar \cX),\ \text{and} \\
 \fR \to \hat   K_0 (MR(\G_k,\QQ_p)),\  \cX \mapsto  \chi_p (\bar \cX).
\end{align*}

\subsection{Cylinders and motivic measure}

Let $\cX$ be a smooth DM stack $\cX$ of finite type and pure dimension $d$.

\begin{defn}
Let $n \in \ZZ_{ \ge 0}$.
A subset $A \subset |\cJ _ \infty \cX|$ is said to be an $n$-\textit{cylinder} if 
$A = \pi _n ^{-1} \pi _n (A)$ and $\pi_n (A)$ is a constructible subset. 
A subset $A \subset |\cJ _ \infty \cX|$ is said to be a \textit{cylinder} if 
it is an $n$-cylinder for some $n$.
\end{defn}

The collection of cylinders is stable under finite unions and finite intersections. 
For a cylinder $A$ and $n \in \ZZ_{\ge 0}$ such that $A$ is an $n$-cylinder,
we define 
\[
\mu_\cX (A):=\{\pi _n (A)\} \LL^{-nd} \in \fR.
\]
This is independent of $n$, thanks to Lemma \ref{lem-trancationpiecewise}.

\begin{rem}\label{rem-difference-Ld}
In some articles, $\tilde \mu_\cX (A)$ is defined to be $\{\pi _n (A)\} \LL^{-(n +1) d}$
instead of $\{\pi _n (A)\} \LL^{-n d}$.
However the difference is just superficial.
\end{rem}

\begin{lem}\label{lem-trancationpiecewise}
Let $\cX$ be a smooth DM stack of pure dimension $d$.
Let  $p :\cJ_{n+1} \cX \to  \cJ_n \cX$ be the natural projection.
Then for every geometric point $z\in (\cJ_n \cX)(K)$, the fiber $p^{-1}(z)$ is isomorphic to 
 $\AA^d_{K}$.
\end{lem}

\begin{proof}
We may assume that $\cX$ is a quotient stack $[M/G]$.
Then $\cJ_n ^l\cX \cong \coprod_{a \in \Conj^l(G)} [J_{nl}^{(a)} M / C_a ]$.
Therefore it suffices to show that the fiber of the morphism $J_{(n+1)l}^{(a)} M \to
J_{nl}^{(a)} M$ over any geometric point 
is isomorphic to an affine space of dimension $d$.
Now we can take a completion of 
$M$ at a $K$-point $w$, $\hat M= \Spec K[[x_1 ,\dots, x_d]]$.
Let $a:\bmu_l \hookrightarrow G$ be an embedding.
By the natural morphism $J_{nl}^{(a)} \hat M \to J_{nl}^{(a)} M$,
the space $(J_{nl} ^{(a)}\hat M)(K)$ is identified with the fiber of $J_{nl}^{(a)}M \to M$ over $w$. 
 Assume that the $\bmu_l$-action on $M$ through $a$ is diagonal and that $\zeta \in \bmu_l$ sends $x_i$ to $\zeta^{a_i}x_i$ with $1 \le a_i \le l$.
Then $(J_{nl}^{(a)} M)(K)$ parameterizes
the homomorphisms $ K[[x_1 ,\dots, x_d]] \to K[[t]]/t^{nl+1} $ which send $x_i$ to an element of the form
\[
 c_0 t^{a_i} +   c_1 t^{l +a_i} +  c_2 t^{2l + a_i} \cdots + c_{n-1} t^{(n-1)l+a_i}, \ c_i \in K.
\] 
Therefore the fiber of $J_{(n+1)l}^{(a)} M \to
J_{nl}^{(a)} M$ over every geometric point
is isomorphic to an affine space of dimension $d$.
\end{proof}

It is obvious that $ \mu_\cX$ is a finite additive measure:

\begin{prop}
 If $A$ and $A_i$ $(i=1,\dots,n)$ are cylinders such that $A = \coprod_{1 \le i \le n} A_i$, then 
\[
 \mu_\cX (A) = \sum _{i=1}^n \mu_\cX (A_i) .
\]
\end{prop}

\subsection{Integrals of measurable functions}\label{subsec-integrals-of-functions}

\begin{defn}
Let $A \subset |\cJ_\infty \cX|$ be a subset. 
A function $F : A \to \fR$ is said to be {\em measurable}
if there are countably many cylinders $A_i$ such that $A = \coprod _i A_i$
and the restriction of $F$ to each $A_i$ is constant. 
\end{defn}

We define the integral of a measurable function $F: A \to  \fR $ as follows;
\[
 \int _A F d \mu _\cX := \sum  F(A_i) \cdot \mu_\cX (A_i)  \in \bar \fR.
\]

\begin{defn}
Let $A \subset |\cJ_\infty \cX|$ be a cylinder. 
If $n \in \ZZ_{\ge 0}$ is such that 
$A$ is an $n$-cylinder, 
we define 
\[
\codim A := \codim ( \pi _n(A), |\cJ_n \cX|) .
\]
\end{defn}

\begin{lem}
The integral of a measurable function $F$ is independent of the choice of $A_i$.
\end{lem}

\begin{proof}
Let $A_i$, $B_i$, $i \in \NN$ be cylinders such that $F|_{A_i }$ and $F|_{B_i}$ are constant
for all $i$ and $A = \coprod _{i \in \NN} A_i = \coprod _{i \in \NN}B_i$.
We assume that $\{B_i\}$ is a refinement of $\{A_i\}$.
The general case can be reduced to this case.
Let $S_i := \{ j\in \NN | B_j \subset A_i \}$,
 let $\eta _1 , \dots , \eta _u $ be the generic points of
the irreducible components of $A_i$ and let $ j_ v \in S_i$, $1 \le v \le u$ be
such that $\eta _v \in B_{j_v} $. 
Then we see that for $j' \in S_i \setminus \{j_1, \dots, j_u\}$,
\[
 \codim B_{j'} \ge  \codim (A_i \setminus \bigcup_{1 \le v \le  u} B_{j_v} ) > \codim A_i .
\]
Repeating this argument, we see that for every $m \in \ZZ$, there are at most finitely many $j \in S_i$ such that
$ \codim B_j <m $.
It means that for every $n \in \ZZ$, there is a finite subset $S_{i,n} \subset S_i$
such that for $j\in S_i \setminus S_{i,n}$, $\dim \mu_\cX (B_j) < n$.
For every $n$, we have
\[
 \mu_\cX (A_i) \sim _n \mu_\cX (\coprod _{ j \in S_{i,n} } B_j) \sim _n
 \sum _{ j \in S_{i} }  \mu_\cX (B_j) .
\]
Hence $ \mu_\cX (A_i) = \sum _{ j \in S_{i} }\mu_\cX ( B_j)  $.
This proves the lemma. (In fact, $S_i$ is a finite set. See \cite[Lemma 2.3]{looi})
\end{proof}

Below, we consider such functions as $\LL^{\ord \cI _{\cZ}}$. Here $\ord \cI _{\cZ}$ is the 
order function associated to the ideal sheaf $\cI_{\cZ}$ of a closed substack $\cZ \subset \cX$
and takes values in $\ZZ_{\ge 0}$ outside $|\cJ_\infty \cZ| \subset |\cJ_\infty \cX|$ and
the infinity on $|\cJ_\infty \cZ|$. Thus a function $\LL^{\ord \cI _{\cZ}}$ is not defined 
on $|\cJ_\infty \cZ|$ (at least, as a $\fR$-valued function). However, 
if $\cZ \subsetneqq \cX$, then
this does not cause any problem, because $|\cJ_\infty \cZ|$ is too small to affect integrals.

\begin{defn}
A subset $A \subset |\cJ_\infty \cX|$ is said to be {\em negligible} if 
there are cylinders $A_i$, $i \in \NN$ such that $A =  \bigcap _i A_i$
and $\lim _{i \to \infty}\codim A_i = \infty$. 
\end{defn}

\begin{prop}\label{proposition-measure-zero}
Let $\cZ \subset \cX$ be a locally closed substack of
dimension $d' < d$. Then a subset $|\cJ _\infty \cZ| \subset |\cJ_\infty  \cX|$ is 
negligible.
\end{prop}

\begin{proof}
We have
\[
 |\cJ_\infty \cZ| = \bigcap_{n \ge 0} \pi_n ^{-1} \pi_n (|\cJ_\infty \cZ|) .
\]
From Lemma \ref{lemma-Greenberg},
$\pi_n ^{-1} \pi_n (|\cJ_\infty \cZ|)$ are cylinders, and from Lemma \ref{substack-trancated-arc-fiber}, 
\[
 \lim _{n \to \infty} \codim \pi_n ^{-1} \pi_n (|\cJ_\infty \cZ|) = \infty .
\]
Hence $|\cJ_\infty \cZ|$ is negligible.
\end{proof}

Removing a negligible subset from the domain or adding one to the domain does 
not change the value of  integrals:

\begin{prop}\label{proposition-negligible-is-negligible}
Let $F :A \to \fR $ be a function and $B\subset A$ a negligible subset.
If $F$ is measurable, then $F |_{A \setminus B}$ is measurable and 
\[
 \int _A F d \mu_\cX = \int _{A \setminus B} F d \mu_\cX .
\]
\end{prop}

\begin{proof}
Let $A_i$, $i \in \NN$ be cylinders such that $A = \coprod _i A_i$ and $F|_{A_i}$
are constant for all $i$. Let $B_i$, $i \in \NN$ be cylinders such that $B = \bigcap _i B_i$ and $\lim _{i \to \infty} \codim B_i \to  \infty$. 
We may assume that $B_1  = A$ and  $B_{i+1} \subset B_i$ for all $i$.
We put $C_i : = B_i \setminus B_{i+1}$. Then $C_i$ are mutually disjoint cylinders such that
$\coprod _{i\in \NN} C_i = A \setminus B$ and $\lim_{i \to \infty} \codim C_i = \infty$.
For $i,j\in \NN$, we put $A_{ij} := A_i \cap C_j$.
These are cylinders with $\coprod _{i,j} A_{ij}= A \setminus B$.
Since $F |_{A_{ij}}$ are constant for all $i,j$, $F |_{A \setminus B}$ is measurable.

We have
\begin{align*}
\int _{A \setminus B}F d \mu_\cX = \sum _i \sum _j F(A_{ij}) \cdot \mu_\cX (A_{ij}).
\end{align*}
Therefore, to prove the equation of the proposition, it suffices to show
$\mu_\cX (A _i) = \sum_j \mu _\cX (A_{ij})$.
By definition, we have $A_i \setminus B_{j_0+1} = \coprod _{j \le j_0} A_{ij}$.
For every $n \in \ZZ$, if $j_0 $ is sufficiently large, then
\begin{align*}
\mu_\cX (A _i) \sim _n  \mu_\cX (A_i \setminus B_{j_0+1}) = \mu_\cX ( \coprod _{j \le j_0} A_{ij} ) 
\sim_n \sum_j \mu _\cX (A_{ij}) .
\end{align*}
Hence $\mu_\cX (A _i) = \sum_j \mu _\cX (A_{ij})$.
\end{proof}

By abuse of terminology,
we say that $F$ is a measurable function on a subset $A$ even if
$F$ is a measurable function defined only on $A \setminus B$ with $B$ negligible, and write 
\[
\int_{A} F d \mu_\cX = \int_{A \setminus B} F d \mu_\cX.
\]

\begin{rem}
The definition of the measurable function in this paper differs from that in \cite[Appendix]{DL-quotient}.
\end{rem}

\subsubsection{Lemmas}

\begin{lem}\label{lemma-Greenberg}
Let $\cZ$ be a DM stack of finite type. Then
there is a monotone increasing function  $\phi: \ZZ_{ \ge 0 } \to \ZZ_{ \ge 0 }$ such that $\phi (n) \ge n$ for all $n$ and 
\[
 \pi _n (|\cJ_\infty \cZ |)= \Im(|\cJ_{\phi(n)} \cZ| \to |\cJ_n \cZ| ) .
\]
In particular, $\pi _n (|\cJ_\infty \cZ |) \subset |\cJ_n\cZ|$ is a constructible subset.
\end{lem}

In the case where $\cZ$ is a scheme, this was proved by Greenberg \cite{Greenberg-henselian}
by using Newton-Hensel lemma (called Newton's lemma in \cite{Greenberg-henselian}).
We prove Lemma \ref{lemma-Greenberg} by using the equivariant version of 
Newton-Hensel lemma (Lemma \ref{lemma-generalized-Hensel}).

\begin{proof}
The second assertion of the lemma follows from the first and 
Chevalley's theorem for stacks \cite[Th\'eor\`eme 5.9.4]{LMB}.

We may assume that $\cZ \cong [Z/G]$ with $Z=\Spec R$ an affine scheme and $G$ a finite group.
Furthermore we may assume that $k$ contains all $l$-th roots of unity for $l$ prime to the characteristic of $k$ such that there is
an element $g \in G$ of order $l$. 
From  Proposition \ref{lem-verifyDMstack}, the first assertion of the lemma is equivalent to the following:
\begin{itemize}
\item[$\bigstar $] For every $a:\bmu_l \hookrightarrow G$ and for every $ 0 \le n < \infty $, there is
a monotone increasing function $\phi: \ZZ_{\ge 0} \to \ZZ_{\ge 0}$ such that  for every $n$,
  $\phi (n) \ge n$ and
\[
 \Im (J^{(a)} _{\phi(n)l} Z \to J^{(a)}_{nl} Z) = \Im (J^{(a)} _{\infty} Z \to J^{(a)}_{nl} Z).
\]
\end{itemize}

We prove $\bigstar $ by the induction on $\dim Z$. 
There exists a $\bmu_l$-equivariant embedding of $Z$ into $\AA^ d = \Spec k[x_1 ,\dots,x_d] $
 on which $\bmu_l$ acts diagonally. 
 Let $I   \subset k[x_1 ,\dots, x_d]$ be the defining ideal of $Z$.
 For some positive integer $s$, we have $\sqrt{I} \supset I \supset (\sqrt{I})^s $. Let $Z _{\mathrm{red}}$ and $Z_{\mathrm{red}}^s$ be the closed subschemes defined
by $\sqrt{I}$ and $ (\sqrt{I})^s $ respectively. 
Then for every $m$, we have
\[
 J^{(a)}_{ml} Z_{\mathrm{red}} \subset  J^{(a)}_{ml} Z \subset J^{(a)}_{ml} Z_{\mathrm{red}} ^s \subset J^{(a)}_{ml} \AA^ d
\]
and 
\[
 J_\infty^{(a)} Z_{\mathrm{red}} = J_\infty^{(a)} Z = J_\infty ^{(a)} Z_{\mathrm{red}}^s.
 \]
Hence it suffices to show $\bigstar$ for $Z _{\mathrm{red}}$ and $Z_{\mathrm{red}}^s$.
If $m' $ is such that $m'l + 1 \ge mls +s$, then
we have $ \Im (J_{m'l}^{(a)}Z_{\mathrm{red}}^s \to J_{ml}^{(a)}\AA^d)  \subset J_{ml}^{(a)} Z_{\mathrm{red}}$.
Hence it suffices to show $\bigstar$ only for $Z$ reduced. 

Then  $\bigstar$ clearly holds if $Z$ is of dimension zero.

We can also assume that $Z$ is irreducible: Let $Z_1, \dots, Z_q$ be the irreducible components of $Z$
and let $W:= \bigcup _{i\ne j} Z_i \cap Z_j$. Since $\dim W < \dim Z$, by the inductive assumption,
$\bigstar$ holds for $W$. Every $\gamma \in J_\infty Z \setminus J_\infty W$ lies in only one component $Z_i$.
Thus Assertion $\bigstar$ for $Z$ follows from $\bigstar$  for $Z_1 ,\dots, Z_q$.

Suppose that $Z$ is irreducible, reduced  and of positive dimension.
Let $f_1, \dots, f_s \in k[x_1 ,\dots ,x_d]^{\bmu_l}$  be $\bmu_l$-invariant 
polynomials defining  $Z$ and
 $r := \codim (Z ,\AA^d) \le s$. 
Reordering suitably $f_1, \dots,f_s$,
we have that $Z$ is an irreducible component of $V:=V(f_1 ,\dots ,f_r)$.
Let
\[
J := \det \left(\frac{\partial f_i}{\partial x_j}\right)_{1\le i \le r,\ 1 \le j \le r}.
\]
Since $k$ is perfect and $Z$ is reduced, $Z$ is generically smooth.
Therefore, if necessary, reordering variables $x_1,\dots,x_d$, we may assume 
that the subscheme $S:=V(f_1 ,\dots, f_s , J) \subset Z$  is of dimension $< \dim Z$:
Let  $W$ be the intersection of $Z$ and the closure of $V \setminus Z$ and let $X \subset Z$ be
the scheme-theoretic union of $W$ and $S$, which is of dimension $< \dim Z$.
By the inductive assumption, 
there exists a monotone increasing function $\psi$ such that
\[
 \Im(J_{\infty }^{(a)} X \to  J_{n l }^{(a)}X  ) = \Im(J_{\psi(n) l }^{(a)} X \to  J_{n l }^{(a)} X ).
\]

Let $K \supset k$ be a field extension.
A $K$-point $\gamma$ of $J_\infty ^{(a)} \AA^d$ corresponds to a 
$\bmu_l$-equivariant $K[[t]]$-algebra homomorphism 
\[
 \gamma^* : K[[t]][x_1,\dots ,x_d] \to K[[t]].
 \]
Suppose that $\zeta \in \bmu_l$ sends $x_i$ to $\zeta ^{a_i} x_i $, $0 \le  a_i \le l-1$.
Let $\gamma \in (J_\infty ^{(a)} \AA^d)(K)$ be such that $\pi_{nl} (\gamma) \in J_{nl}^{(a)} Z \setminus 
J_{nl}^{(a)} X $ and let $nl\le  e < (n+1)l $ be the unique integer such that 
$e + \sum _{i = 1}^r a_i \equiv 0 \mod l $. 
Since $\pi_{nl} (\gamma) \notin J_{nl} ^{(a)}S$,
 we have 
\[
\gamma^*(J) \not \equiv 0 \mod (t^{nl+1})\ (\text{and }\mod (t^{e+1})) .
\]
From  Lemma \ref{lemma-generalized-Hensel}, if  $\pi_{ml}(\gamma) \in J_{ml}^{(a)} V$ with $ml+1 \ge 2 (e + \sum _{i = 1}^r a_i) +l$, then there is $\gamma' \in J_\infty ^{(a)} V $ such that
 $\pi_{nl} (\gamma')= \pi_{nl} (\gamma)$.
Moreover $\pi_{ml}(\gamma)$ must lie in $ J_{ml}^{(a)} Z$ and $\gamma'$ in $J_\infty ^{(a)}Z$.
Hence there exists a monotone increasing function $\tau $ such that
\[
 \Im (J_\infty ^{(a)} Z \to J_{nl}^{(a)} Z )  \setminus J_{nl}^{(a)}X=\Im (J_{\tau(n)l}^{(a)} Z \to J_{nl}^{(a)} Z ) \setminus J_{nl}^{(a)}X.
\]

We show that $\bigstar$ holds for $\phi := \tau \circ \psi $.
Let $ v \in  \Im ( J_ {\phi (n)l}^{(a)}Z \to J_{nl}^{(a) } Z ) $ be an arbitrary point. 
We have to show that $v \in \Im (J_\infty ^{(a)}Z \to J_{nl}^{(a) } Z )$.
This holds, either if $v \notin  J_{nl}^{(a) } X$ or if 
$ v \in \Im ( J_{\psi (n)l} ^{(a)}X \to J_{nl}^{(a) } X  ) $.
If it is not the case, there exists $w \in \Im ( J_ {\phi (n)l}^{(a)}Z \to J_{\psi (n)l}^{(a) } Z ) \setminus J_{\psi (n)l}^{(a)}X$ which maps to $v$.
By the definition of $\tau$, we have $w \in \Im (J_\infty ^{(a)}Z \to J_{\psi(n)l}^{(a)}Z)$ and 
$v \in \Im (J_\infty ^{(a)}Z \to J_{nl}^{(a)}Z)$.
\end{proof}

\begin{lem}[Equivariant Newton-Hensel lemma]\label{lemma-generalized-Hensel}
Let $l$ be a positive integer prime to the characteristic of $k$. 
Suppose that $k$ contains all $l$-th roots of unity.
Suppose that $\bmu_l$ acts on $k$-algebras $k[[t]][x_1 , \dots, x_d]$ and  $k[[t]]$ by
\begin{align*}
\bmu _l \ni \zeta :& t \mapsto \zeta t \\
& x_i \mapsto  \zeta ^{a_i} x_i \ (a_i \in \{0,1,\dots, l-1\}) \text{ and } \\
  \bmu _l \ni \zeta :& t \mapsto \zeta t \text{ respectively}.
\end{align*}
Let  $\Lambda$  be the set of $\bmu_l$-equivariant $k[[t]]$-algebra homomorphisms
\[
k[[t]][x_1,\dots ,x_d] \to k[[t]].
\]
Let $f_i$ $(i=1,2,\dots,r \le d)$ be elements in the invariant subring $k[[t]][x_1,\dots,x_d]^{\bmu_l}$
and  let 
\[
J : = \det \left(\frac{\partial f_i}{\partial x_j}\right)_{1\le i \le r,\ 1 \le j \le r}\in k[[t]][x_1,\dots,x_d].
\] 
Suppose that for a positive integer $e$ such that $e':=e + \sum_{j=1}^r a_j$ is divisible by $l$
  and for $\eta \in \Lambda$,
we have 
\[
\eta (J) \not \equiv 0 \mod (t^{e+ 1})
\]
and for $i=1,2,\dots,r$,
\[
 \eta (f_i) \equiv 0 \mod (t^{2 e' + l})  .
\]
Then there is $\theta \in \Lambda$ such that for every $i$,
\[
\theta (f_i) =0 \text{ and } \theta (x_i) \equiv \eta (x_i) \mod ( t^{e + l}) .
\]
\end{lem}

\begin{proof}
Let $\gamma : k[[t]][x_1 ,\dots,x_d]\to k[[t]]$ be a $k[[t]]$-algebra homomorphism.
Then $\gamma \in \Lambda$  if and only if for every $i$, $\gamma(x_i) \in  t^{a_i } \cdot k[[t^l]]$.
Consider an $k((t))$-automorphism
\[
 \alpha : k((t))[x_1,\dots,x_d] \to k((t))[x_1,\dots,x_d], \ x_i \mapsto t^{-a_i}x_i .
\]
If $\tilde \gamma :k((t))[x_1,\dots,x_d] \to k((t))$ is a $k((t))$-algebra homomorphism which is the extension of $\gamma \in \Lambda$,
then for every $i$, $\tilde \gamma \circ \alpha(x_i) $ lies in a subring $k[[t^l]] \subset k[[t]]$ and
hence $\tilde \gamma\circ \alpha (k[[t^l]][x_1 ,\dots, x_d]) \subset k[[t^l]]$.
We define a $k[[t^l]]$-algebra homomorphism 
\[
\gamma' : k[[t^l]][x_1,\dots,x_d] \to k[[t^l]]
\]
to be the restriction of $\tilde \gamma \circ \alpha$.
Let $\Lambda'$ be the set of $k[[t^l]]$-algebra homomorphisms $ k[[t^l]][x_1,\dots,x_d] \to k[[t^l]]$. Then the map 
\[
\Lambda \to \Lambda', \ \gamma \mapsto \gamma'
\]
is  bijective.

The invariant ring $k[[t]][x_1,\dots,x_d]^{\bmu_l}$ is generated as a $k[[t^l]]$-algebra  by the monomials
$t^b x_1^{b_1} \dots x_d^{b_d}$ such that $b + \sum_{i=1}^d a_i b_i \equiv 0 \mod l$.
Since 
\begin{align*}
\alpha^{-1}(t^b x_1^{b_1}\cdots x_d ^{b_d}) 
= t^{b + \sum a_ib_i }x_1^{b_1} \cdots x_d^{b_d} \in k[[t^l]][x_1 , \dots ,x_d],
\end{align*}
$f'_i:=\alpha^{-1}(f_i)$ lie in $k[[t^l]][x_1,\dots,x_d]$.
Let 
\[
J':=\det \left(\frac{\partial f'_i}{\partial x_j}\right)_{1\le i \le r,\ 1 \le j \le r}.
\]
For any $g \in k((t))[x_1,\dots,x_d] $, by easy calculation, we can see
\[
 \frac{ \partial \alpha^{-1}(g) }{\partial x_i} = t^{a_i} \alpha^{-1} \left(\frac{ \partial g }{\partial x_i} \right).
\]
Hence we have
\[
  J' =  t ^{\sum_{j=1}^r a_j}\alpha^{-1} ( J) .
\]
We can assume that $ (\eta ' (J')) = (t^{e'})  $, or equivalently that
$(\eta (J)) = (t ^e)$: For if $\eta'(J')=(t^c) \subset k[[t^l]]$,  we replace
$f_1$ and $f'_1$ with $t^{e'-c} f_1$ and  $t^{e'-c}f'_1$.

We have 
\[
 \eta'(f'_i) = \eta (f_i) \in  (t^{2e' + l}) \subset k[[t^l]].
\]
Now from the (non-equivariant) Newton-Hensel lemma \cite[Chapitre III, \S 4, Corollaire 3]{Bourbaki_Algebre_Commutative}, there exists $\theta' \in \Lambda'$ such that
\[
 \theta ' (f ' _i) =0 \text{ and } \theta'(x_i) \equiv \eta'(x_i)  \mod ( t^{e' + l}) .
\]
Let $\theta$ be the element mapping to $\theta'$ by the bijection above $\Lambda \to \Lambda'$.
Then $\theta (f_i)=0$ for every $i$. Moreover since $\theta (x_i) = \theta'(x_i) t^{a_i}$ and
$\eta(x_i) = \eta'(x_i)t^{a_i}$, we have that for every $i$,
\[
 \theta (x_i) \equiv \eta(x_i) \mod (t^{e' + l+a_i }) ,
\]
hence 
\[
 \theta (x_i) \equiv \eta(x_i) \mod (t^{e + l}) .
\]
\end{proof}

\begin{lem}\label{substack-trancated-arc-fiber}
Let $\cZ$ be a DM stack of finite type and dimension $d'$.
\begin{enumerate}
\item  For every $0 \le n < \infty$,
every fiber of $\pi _{n+1} (\cJ_\infty \cZ) \to \pi_n (\cJ_\infty \cZ) $ is of dimension $\le d'$.
\item  $ \dim \pi _n (\cJ_\infty \cZ) \le d' (n+1) $.
\end{enumerate}
\end{lem}

\begin{proof}
The second assertion is a direct consequence of the first. 
Now we may assume that $k$ is algebraically closed and $ \cZ$ is a quotient stack $[Z/G]$.
From Proposition \ref{lem-verifyDMstack}, it suffices to show that for every embedding 
$a : \bmu_l \hookrightarrow G$, 
every closed fiber of a morphism
\begin{equation}\label{equation-JainftyZ}
 \pi_{(n+1)l}(J_{\infty}^{(a)}Z) \to \pi_{nl}(J_{\infty}^{(a)}Z)
\end{equation}
is of dimension $\le d'$. 
Take a $\bmu_l $-equivariant embedding 
\[
Z \hookrightarrow \AA^d= \Spec k[x_1,\dots,x_d]
\]
 where 
$\bmu_l$ acts on $Z$ through $a$ and on $\AA^d$ by 
\[
\bmu_l \ni \zeta : x_i \mapsto \zeta ^{a_i}x_i , \ 0 \le a_i \le l-1.
\]
Let $ \mathbf{f}=(f_1,\dots, f_r)$ be a system of polynomials in $k[x_1,\dots,x_d]$ which defines
$Z$. Then $(J_\infty ^{(a)} Z)(k)$
is identified with 
\[
 \{(\phi _1 , \dots ,\phi _d) | \phi _i \in  t ^ {a_i}\cdot k[[t^l]] \subset k[[t]], \ 
 \mathbf{f}(\phi _1 , \dots ,\phi _d)=0 \} .
\]
Let $\gamma \in J_\infty^{(a)} Z$ correspond to $(\phi _1 , \dots ,\phi _d)$.
Then the fiber of (\ref{equation-JainftyZ}) over $\pi_{nl} (\gamma)$ is
identified with 
\[
B:= \{ (\bar \psi_1 ,\dots ,\bar \psi _d )| \mathbf{f}(\dots, \phi_i + t ^{a_i + nl} \psi _i , \dots) = 0, \, 
  \psi _i \in k[[t^l]]   \} .
\]
Here $\bar \psi$ is the image of $\psi$ by $k[[t]] \twoheadrightarrow  k=k[[t]]/t$. This is contained in 
\[
 B':=\{ (\bar \psi_1 ,\dots ,\bar \psi _d )| \mathbf{f}(\dots, \phi_i + t ^{a_i + nl} \psi _i , \dots) = 0, \, 
  \psi _i \in k[[t]]   \} .
\]
The equations $ \mathbf{f}(\dots, \phi_i + t ^{a_i + nl} \psi _i , \dots) = 0 $ define a closed subscheme 
\[
Y \subset \Spec k[[t]][[\psi _1 ,\dots,\psi_d  ]].
\]
The generic fiber of the projection $Y \to \Spec k[[t]]$ is isomorphic to $Z \otimes _k k((t))$.
Then $B'$ is contained in the intersection of the special fiber and the closure of the generic fiber.
Hence $\dim B \le \dim B' \le d'$. 
\end{proof}

\subsection{Motivic integration over  singular varieties}
We review the motivic integration over singular varieties which was studied by Denef and Loeser \cite{germs}.
They assumed that the base field is of characteristic zero.
Their arguments however apply to an arbitrary perfect
field, as verified by Sebag \cite{Sebag} in a more general situation. 
Although they considered only schemes, we can simply generalize the theory to algebraic spaces.

Comparing a DM stack and its coarse moduli space is an interesting problem. Even if the stack is 
smooth, the coarse moduli space is not generally smooth. Therefore, we consider not only the motivic integration
over smooth DM stacks, but also that over singular varieties.

Let $X$ be a reduced variety of pure dimension $d$.
If $X$ is not smooth, then fibers of $J_{n+1} X \to J_{n} X$ are not generally isomorphic to $\AA^d$.
From Greenberg's theorem \cite{Greenberg-henselian}, $\pi _n (J_\infty X)$ is a constructible subset. 
Denef and Loeser proved that
every fiber of $\pi _{n+1} (J_\infty X) \to \pi _n(J_n X)$ is of dimension $\le d$ \cite[Lemma 4.3]{germs}.

For an ideal sheaf $\cI \subset \cO_X$ and $\gamma \in (J_\infty X)(K)$ with $K$ a field, 
we define the order of $\cI$ along $\gamma$ to be $\ord \cI (\gamma):= n$ if
$\gamma^{-1} \cI = (t^n) \subset K[[t]]$. By convention, we put $\ord \cI (\gamma):= \infty$ if
$\gamma^{-1} \cI = (0)$. Thus we have the order function associated to $\cI$,
\[
 \ord \cI : J_\infty X \to \ZZ _{\ge 0} \cup \{\infty\} .
\]
Let $\Jac_X$ be the Jacobian ideal sheaf of $X$, that is, the $d$-th Fitting ideal of 
$\Omega _{X/k}$. 
We define
\[
J_n ^ \diamondsuit   X   := \{ \pi_n (\gamma) | \gamma \in J_\infty X ,\ \ord \Jac _X (\gamma) <n \} 
\subset \pi _n (J_\infty X).
\]
This is a constructible subset of $J_n X$.
The fibers of
 $\pi _{n+1} (J_\infty X) \to \pi _n (J_\infty X)$ over $J_n ^ \diamondsuit   X $ are isomorphic to a $d$-dimensional affine space
 (see \cite[Lemma 9.1]{looi}).

\begin{defn}
A subset $A \subset J_\infty X$ is called an {\em $n$-cylinder} if $A = \pi _n ^{-1}\pi _n (A)$
and $\pi _n(A)$ is a constructible subset in $ J_n ^\diamondsuit X$. A subset is called a {\em cylinder} if it is an $n$-cylinder
for some $n \in \ZZ _{\ge 0}$.
For an $n$-cylinder $A$, we define $\codim A :=   (n+1)d - \dim \pi_n(A)$.
\end{defn}

For an $n$-cylinder $A$, 
we define 
\[
\mu_X(A):= \{ \pi_n(A)\}\LL^{-nd} \in \fR.
\]
As in the case of smooth stacks, 
we say that a function $F: J_\infty X \supset A \to \fR$ is {\em measurable} if there are 
countably many cylinders $A_i$ such that $A = \coprod _{i } A_i$ and
the restriction of $F$ to each $A_i$ is constant. 
For a measurable function $F$, we define
\[
 \int _A F d \mu_\cX := \sum  F(A_i)\cdot \mu_\cX (A_i) \in \bar \fR.
\]

\begin{defn}\label{defn-negligible-variety}
A subset $A \subset J_\infty X$ is said to be a {\em negligible} if there are constructible subsets
$C_n \subset \pi _n (J_\infty X)$, $n \in \NN$ such that
$A = \bigcap _{n \in \NN} \pi _n ^{-1} (C_n)$ and $ \lim _{n \to \infty} \dim C_n - dn = - \infty$.
\end{defn}

When $X$ is smooth, this definition coincides with that in the preceding subsection.
For a subvariety $Z \subset X$ of positive codimension, 
$J_\infty Z $ is a negligible subset of $J_\infty X$.
We indeed have 
\[
 J_\infty Z =  J_\infty Z _{\mathrm{red}} =\bigcap_{n \in \NN} \pi _n^{-1} \pi _n (J_\infty Z_{\mathrm{red}}), 
\]
where $Z_{\mathrm{red}} \subset X$ is the reduced subscheme associated to $Z$. 
As Denef and Loeser proved, the constructible subset $\pi _n (J_\infty Z_{\mathrm{red}})$ is of dimension
$\le (n+1) \dim Z$.

\begin{prop}\label{proposition-negligible-is-negligible(variety)}
Let $F :A \to \fR$ be a function and $B\subset A$ a negligible subset.
Then if $F$ is measurable, then $F |_{A \setminus B}$ is measurable and
\[
 \int _A F d \mu_\cX = \int _{A \setminus B} F d \mu_\cX .
\]
\end{prop}

\begin{proof}
Let $E_i := \{ \gamma \in |J_\infty X| | \ord \Jac _X (\gamma) = i\}$, $i \ge 0$ and
let $C_n\subset \pi_n (J_\infty X)$ $n \in \NN$ be constructible subsets such that
$ B = \bigcap _{n} \pi_n^{-1} (C_n) $ and $\lim_{n \to \infty} \dim C_n -dn =- \infty $.
Replacing $A$ and $B$ with $A \cap E_i$ and $B \cap E_i$,
we may assume that $A, B \subset E_i$ for some $i$.
Subsets  $C'_n := \pi_n^{-1} (C_n) \cap E_i$ are cylinders 
and $\lim _{n \to \infty} \codim C'_n=\infty $.
 Now we can prove the assertions 
by the same argument as in the proof of Proposition \ref{proposition-negligible-is-negligible}.
\end{proof}

Again, by abuse of terminology, we say that $F$ is a measurable function on $A$ even if
$F$ is a measurable function defined only on $A \setminus B$ with $B$ negligible.

Let $f :Y \to X$ be a proper birational morphism of reduced varieties of pure dimension 
and let
 $X' \subset X$ and $Y ' \subset Y$ be  proper closed subsets such that $f : Y \setminus Y' \cong X \setminus X'$.
Then from the valuative criterion for the properness, 
the map $f_\infty:J_\infty Y \setminus J_\infty Y' \to J_\infty X \setminus J_\infty X'$ is bijective.
In other words, the map 
$f_\infty:J_\infty Y \to J_\infty X $ is bijective outside negligible subsets.
The most fundamental theorem in the theory is 
the following transformation rule (the change of variables formula) which
 describes the relation of $\mu_X$ and $\mu_Y$. 
This was by Kontsevich \cite{Orsay}, Denef and Loeser \cite{germs}, and Sebag \cite{Sebag}.

\begin{thm}\label{TransformationDenefLoeser}
Let $f : Y \to X$ be a proper birational morphism of reduced varieties of pure dimension.
Assume that $Y$ is smooth. Let $\Jac _f \subset \cO_Y$ be the Jacobian ideal of the morphism $f$, that is, the $0$-th 
Fitting ideal of $\Omega _{Y/X}$.
Let $F: J_\infty X  \supset A\to  \fR$ be a measurable function. 
Then $F$ is measurable if and only if $(F \circ f_\infty)\cdot \LL ^{ - \ord \Jac_f}$ is measurable. 
If they are measurable, then we have
\[
 \int _{A} F d \mu_X = \int_{f_\infty ^{-1}(A)} (F \circ f_\infty)\cdot \LL ^{ - \ord \Jac_f} d \mu _Y .
\] 
\end{thm}

\begin{proof}[Sketch of the proof]
Let  $\gamma \in (J_\infty X)(K) $.
Suppose that $\gamma$ sends the generic point of $\Spec K[[t]]$ into the locus where 
$f^{-1}$ is an isomorphism.
The theorem is essentially a consequence of the facts that $f_\infty$ is bijective outside negligible subsets and that for $n  \gg 0$, the fiber
$f_n^{-1} (f_n \circ \pi _n (\gamma))$ is isomorphic to an affine space of dimension $\ord \Jac _f (\gamma) $. 
(We generalize this fact to the stack case in Lemma \ref{keylemma}.)
\end{proof}

\subsection{Tame proper birational morphisms and twisted arcs}

In this subsection, we generalize to the stack case the fact that 
for a proper birational morphism $f$ of varieties, $f _\infty$ is bijective outside negligible subsets.

\begin{defn}
A morphism $f:\cY \to \cX$ of DM stacks is said to be {\em tame} if 
for every geometric point $y$ of $\cY$, the kernel of $\Aut (y) \to \Aut(f(y))$ is of order prime to 
the characteristic of $k$. 
\end{defn}

The following are clear.

\begin{lem}
\begin{enumerate}
\item Tame morphisms are  stable under base change. 
\item Every representable morphism is tame.
\item The composite of tame morphisms is tame.
\end{enumerate}
\end{lem}

\begin{defn} 
A morphism $f:\cY \to \cX$ of DM stacks is said to be \textit{birational} if
there are open dense substacks $\cY _0 \subset \cY$ and $\cX_0 \subset \cX$ such that
$f$ induces an isomorphism $\cY_0 \cong \cX_0$.
\end{defn}

For example, 
 given an effective action of a finite group $G$ on an irreducible variety $M$ 
(that is, for $1 \ne g \in G$, $M^g \subsetneqq  M$), then the natural morphism 
from the quotient stack $[M/G]$ to the quotient variety $M/G$ is birational.
More generally, the morphism from a DM stack $\cX$ to its coarse moduli space is
birational if $\cX$ contains an open dense substack which is isomorphic to an algebraic space.

\begin{prop}\label{birational-bijective}
Let $f: \cY \to \cX$ be a tame proper birational morphism of DM stacks.
Let $\cY'  \subset \cY$ and $\cX' \subset \cX$ be closed substacks such that
$f $ induces an isomorphism $\cY \setminus \cY' \cong \cX \setminus \cX'$.
Then the map 
\[
f_\infty: |\cJ _\infty \cY|\setminus |\cJ_\infty \cY'| \to |\cJ _\infty \cX| \setminus |\cJ_\infty \cX'|
\]
is bijective. 
In particular, if $\cY$ and $\cX$ are either a smooth DM stack or a reduced variety of pure dimension, then $f_\infty$ is bijective outside negligible subsets.
\end{prop}

\begin{proof}
A weak version of this lemma was proved in \cite[Lemma 3.17]{Twistedjet}.
A similar argument works in this general case.

Let $K\supset k$ be an algebraically closed field and
 $ \gamma': \cD _{\infty,K} ^{l'} \to \cX$  a twisted arc such that the generic point maps into
$ \cX \setminus \cX' $.  
Let $\cE$ to be the irreducible component of
 $\cD_{\infty,K} ^{l'} \times _{\cX} \cY$ which contains $\Spec K((t)) \times _{\cX} \cY $ and
 $\cD$  the normalization of $\cE$. Then the natural morphism $\cD \to \cY$ is representable.
 The stack $\cD$ is tame and formally smooth over $K$, and the natural morphism
 $\cD \to \cD_{\infty , K}^{l'}$ is proper and birational. Hence $\cD$ must be isomorphic to $\cD_{\infty,K}^l$ for some $l$ which is
 prime to the characteristic of $k$ and a multiple of $l'$.

Choose the isomorphism $\cD \cong \cD_{\infty,K} ^{l}$ so that the morphism
 $\cD_{\infty,K} ^{l} \cong \cD \to \cD_{\infty,K} ^{l'}$ induces the morphism $D_{\infty,K} \to D_{\infty,K}$ of 
 canonical atlases
 defined by $t \mapsto t^{l/l'}$. Thus we have obtained a twisted arc $\gamma : \cD_{\infty,K} ^l \to \cY$.
 In fact, this is the unique twisted arc which maps to $\gamma'$.
 It follows from the universalities of the fiber product and the normalization. 
\end{proof}

\begin{rem}
This proposition does not hold if we consider only {\em non-twisted} arcs $\Spec K[[t]] \to \cX$. 
It is why we have to introduce the notion of twisted jets.
\end{rem}

\subsection{Fractional Tate objects}

We generalize the transformation rule to 
proper tame birational morphisms of DM stacks in \S \ref{subsec-transformation}. 
Then the contribution of automorphisms of points appears in the formula.
It is of the form $\LL^q$ with $q$ a rational number. 
Therefore we extend the ring in which integrals take values so that
it contains fractional powers of $\LL$. 

We have another motivation to consider fractional powers of $\LL$.
In the birational geometry, particularly in the minimal model program, 
we often treat a normal variety $X$ with $\QQ$-Cartier canonical divisor 
(that is, $X$ is $\QQ$-Gorenstein) or more generally a pair $(X,D)$ of a normal variety
and a $\QQ$-divisor such that $K_X +D$ is $\QQ$-Cartier. 
We can define invariants of $X$ or $(X,D)$, integrating a function of the form 
$\LL ^h$ with $h$ a $\QQ$-valued function deriving from $K_X$ or $K_X +D$.
We deal with this subject in the context generalized to DM stacks in the final section.

Replacing a $\ZZ$-valued function $\alpha$ in Definition \ref{def-convergent-stack}
with $\frac{1}{r}\ZZ$-valued function, we obtain a $\frac{1}{r}\ZZ$-convergent stack.
We define the same equivalence relation $\sim$ of the set $(\fR^{1/r})'$ of the isomorphism classes of
$\frac{1}{r}\ZZ$-convergent stacks.
Then we define $\fR^{1/r}$ to be $(\fR^{1/r})'$ modulo $\sim$.
This is also a semiring and endowed with a map
\[
 \dim : \fR^{1/r} \to \frac{1}{r}\ZZ \cup \{-\infty\} .
\]
We denote $\{(\Spec k,1/r)\}$ by $\LL^{1/r}$.
Then we have $(\LL^{1/r})^r= \LL$.
We can naturally consider $\fR^{1/r}$-valued measurable functions $F : A \to \fR^{1/r}$ and
their integrals:
\[
 \int _A F d \mu _\cX = \sum F(A_i) \mu_\cX (A_i) \in \bar \fR^{1/r}:= \fR^{1/r} \cup \{\infty\}.
\]
We similarly define the $\frac{1}{r}\ZZ$-convergent space and the semiring $\fS^{1/r}$ of
equivalence classes of $\frac{1}{r}\ZZ$-convergent spaces.
There is a semiring homomorphism $\fR^{1/r} \to \fS^{1/r}$, $\{\cX\} \mapsto \{\bar \cX\}$.

We simply define a {\em $\frac{1}{r}\ZZ$-indexed Hodge structure} to be a finite dimensional
$\QQ$-vector space $H$ with a decomposition 
\[
H \otimes _\QQ \CC= \bigoplus _{p,q \in \frac{1}{r}\ZZ} H^{p,q}
\]
such that $H^{q,p}=\overline{H^{p,q}}$. 
Then a $\frac{1}{r}\ZZ$-indexed mixed Hodge structure is a finite dimensional $\QQ$-vector space $H$
endowed with a $\frac{1}{r}\ZZ$-indexed weight filtration $W_\bullet $ of $H$ and a $\frac{1}{r}\ZZ$-indexed Hodge
 filtration $F^\bullet$ of $H \otimes _\QQ \CC$.
The associated graded $\bigoplus _{w \in \frac{1}{r}\ZZ} \Gr ^W _ w H$ is a $\frac{1}{r}\ZZ$-indexed Hodge structure.
We denote the category of $\frac{1}{r}\ZZ$-indexed mixed Hodge structures by $MHS^{1/r}$.
For $a \in \frac{1}{r}\ZZ$, the Tate-Hodge 
structure $\QQ(a)$ is defined to be the one-dimensional $\frac{1}{r}\ZZ$-indexed
Hodge structure $H$ such that $H^{-a,-a}$ is the only nonzero component of $H \otimes _\QQ \CC$.
We denote $[\QQ (a)] \in K_0 (MHS^{1/r})$ by $\LL ^{-a}$.
We define also the completion
$ \hat K_0 (MHS^{1/r})$ similarly.
When $k=\CC$, for a $\frac{1}{r}\ZZ$-convergent space $X= (X,\alpha)$,
we define
\[
   \chi_h (X) := \sum _{V \subset X} \chi_h(V)\LL^{\alpha (V)} \in \hat K_0 (MHS^{1/r}).
\]
There are semiring homomorphisms
\begin{align*}
 \fS^{1/r} &\to \hat K_0 (MHS^{1/r}), \ \{X\} \mapsto \chi_h (X) \\
 \fR^{1/r} &\to \hat K_0 (MHS^{1/r}), \ \{\cX\}  \mapsto \chi_h(\bar \cX) . 
\end{align*}

Suppose that $k$ is a finite field and $p$ is a prime number different from the characteristic of $k$.
If there exist $V \in MR (\G_k , \QQ_p)$ such that $V^{\otimes r} \cong \QQ_p (1)$,
then we fix $V$ and denote it by $\QQ_p(1/r)$.
This is pure of weight $-2/r$. 
For a positive integer $a$, we define $\QQ_p (a /r) := \QQ_p (1/r)^ {\otimes a}$ and 
$ \QQ_p (- a /r) $ to be its dual.
We denote an element $[\QQ_p (-1/r)] \in K_0 (MR(\G_k,\QQ_p))$ also by $\LL^{1/r}$.

T. Ito \cite{ito-tetsu} proved that if we replace $k$ with its suitable finite extension, 
then $\QQ_p(1/r)$ exists. He used this to give a new proof of the well-definedness of
stringy Hodge numbers with $p$-adic integrals and the $p$-adic Hodge theory.

When $\LL^{1/r}$ exists, for a $\frac{1}{r}\ZZ$-convergent space $X= (X,\alpha)$,
we define
\[
   \chi_p (X) := \sum _{V \subset X} \chi_p(V)\LL^{\alpha (V)} \in \hat K_0 (MR(\G_k,\QQ_p)).
\]
Then we have semiring homomorphisms
\begin{align*}
 \fS^{1/r} &\to \hat K_0 (MR(\G_k,\QQ_p)), \ \{X\} \mapsto \chi_p (X) \\
 \fR^{1/r} &\to \hat K_0 (MR(\G_k,\QQ_p)), \ \{\cX\} \mapsto \chi_h(\bar \cX) . 
\end{align*}
 
\subsection{Shift number}

Let $\cX$ be a smooth DM stack of pure dimension $d$ and $x \in \cX(K)$ 
a geometric point. Then the automorphism group
$\Aut (x)$ linearly acts on the tangent space $T_x \cX $. 
Let $a:\bmu_l \hookrightarrow \Aut (x)$ be an embedding.
According to the $\bmu_l$-action, $T_x \cX$ decomposes into eigenspaces;
\[
 T_x \cX = \bigoplus _{i =1}^l T_{i,x} .
\] 
Here $T_{i,x}$ is the eigenspace on which $\zeta \in \bmu_l$ acts by the multiplication of 
$\zeta ^i$.
We define
\[
 \sht (a):= d - \frac{1}{l}\sum _{i=1}^l i \cdot \dim T_{i,x} 
 = \frac{1}{l} \sum _{i=1}^l (l-i) \cdot \dim T_{i,x} .
\]
If $a'$ is a conjugacy of $a$, then $\sht (a')=\sht(a)$.
Since the fiber of $|\cJ_0 \cX| \to |\cX|$ over $x$ is identified with 
$\coprod _{\chara (k) \nmid l} \Conj (\bmu_l , \Aut (x))$,
we define $\sht (p) := \sht (a)$ if $p \in |\cJ_0 \cX|$ is the point corresponding to $(x,a)$.
Thus we have a map
\[
\sht :|\cJ_0 \cX| \to \QQ .
\]

Furthermore, $\sht(p)$ depends only on the connected component $\cV \subset \cJ_0 \cX$
in which $p$ lies:
Let $f:S \to  \cV$ be a morphism with $S$ a connected scheme such that
$p : \Spec K \to \cV$ factors as $\Spec K \to S \to \cV$.
Let  $ f' : S \to \cV \to \cX$ be the composite of $f$ and the projection $\cV \to \cX$.
Then the pull-back $(f')^* T\cX$ of the tangent bundle has a $\bmu_{l,k}$-action naturally deriving from
$f$  and decomposes into eigenbudles,
\[
 (f')^{*} T\cX = \bigoplus _{i=1}^l  T_i 
\]
such that $T_{i,x}$ above is a pull-back of $T_i$. 
Then 
\[
 \sht (p) = d - \frac{1}{l}\sum _{i=1}^l i \cdot \rank T_i .
\]
It follows that $\sht (p)$ is constant on $\cV$. 
 We define $\sht (\cV):= \sht (p)$, $p \in |\cV|$.

We denote the composite map $| \cJ_\infty \cX| \to |\cJ_0 \cX| \xrightarrow{\sht} \QQ $ by $\fs _ \cX$.
For a (possibly singular) variety $X$, we denote by $\fs _X $ the constant zero function over $|J_\infty X|$. 
These definitions coincide for smooth varieties. 
In both cases, the function 
\[
\LL^{\fs_\cX} : |\cJ_\infty \cX| \to \fR^{1/r}
\]
 is clearly  measurable.

\subsection{Transformation rule}\label{subsec-transformation}

In this subsection,  we prove the transformation rule generalized to  tame proper birational morphisms.

\begin{defn}
Let $\cX$ be a DM stack and $\cI \subset \cO _{\cX}$ an ideal sheaf.
Let $\cZ \subset \cX$ be the closed substack defined by $\cI$.
We define a function 
\[
\ord \cI :|\cJ_\infty \cX | \setminus |\cJ_\infty \cZ| \to \QQ
\]
 as follows:
Let $\gamma \in |\cJ_\infty \cX| $ and  
$\gamma _K: \cD _{\infty,K} ^l   \to \cX$ its representative. 
Let $ \gamma'_K :D_{\infty,K} \to \cX$ be the composite of $\gamma_K$ and the canonical atlas $ D_{\infty,K}  \to \cD _{\infty,K} ^l $.
Suppose that $( \gamma'_K) ^{-1} \cI = (t ^m) \subset K[[t]]$. 
Then 
\[
 \ord \cI (\gamma) := \frac{m}{l} \in \QQ.
\] 
\end{defn}

Let $\cX$ be a smooth DM stack and $\cI \subset \cO_\cX$ an ideal sheaf such that
the support of $\cO_\cX/\cI$ is of positive codimension. 
Let $r \in \NN$ be such that $\Im (\ord \cI ) \subset \frac{1}{r}\ZZ $.
Then the function 
 \[
 \LL^{\ord \cI} : |\cJ _\infty \cX| \to \fR^{1/r} 
 \]
 is measurable (defined outside a negligible subset $|\cJ_\infty \cZ|$). 

\begin{defn}
Let $f:\cY \to \cX$ be a birational morphism of DM stacks.
We define its {\em Jacobian ideal sheaf} $\Jac _f \subset \cO_\cX$ to be the $0$-th Fitting ideal sheaf of
$\Omega_{\cY/\cX}$.
\end{defn}

\begin{thm}[Transformation rule]\label{thm-transformationstackstack}
Let $\cY$ and $\cX$ be DM stacks of finite type and 
pure dimension $d$ and $f:\cY \to \cX$ a tame proper birational morphism.
Suppose that $\cY$ is smooth and $\cX$ is either a smooth DM stack or a reduced variety.
Let $A \subset |\cJ_\infty \cX|$ be a subset and $F :A \to \fR^{1/r} $ a function
(at least, defined outside a negligible subset).
\begin{enumerate}
\item\label{assertion-measurable} The function $F$ is measurable if and only if $F \circ f_\infty$ is measurable.
\item\label{assertion-formula} Suppose that $F$ is measurable and that the function
$\LL^{ \fs _\cX }|_A$ takes values in $\fR^{1/r}$. Then 
\[
\int _A  F \LL^{ \fs _\cX } d\mu_{\cX} = 
\int _{f_\infty^{-1}(A)} (F \circ f_\infty) \LL^{- \ord \Jac_f + \fs _\cY} d\mu_{\cY}\  \in 
\bar \fR ^{1/r}.
\]
\end{enumerate}
\end{thm}

\begin{proof}
\ref{assertion-measurable}.
Let $\cX'\subset \cX$ be the closed substack over which $f$ is not an isomorphism.
Let $A_i \subset |\cJ_\infty \cX| \setminus |\cJ_\infty \cX'|$, $i \in \NN$ be subsets such that $\coprod _{i\in \NN} A_i = A \setminus |\cJ _\infty \cX'| $ and $F$ is constant over each $A_i$. 
Let $B_i := f^{-1}_\infty (A_i)$. From Lemma \ref{key lemma 2}, $A_i$ is a cylinder if and only if
$B_i$ is a cylinder. Thus $F$ is measurable if and only if $F \circ  f_\infty$ is measurable.

\ref{assertion-formula}.
Let $A_i$ and $B_i$ be as above.
If necessary, taking a refinement of $\{A_i\}_i$,
we may assume that $\fs_\cX$ is constant on every $A_i$ and that
$\fs_\cY$ and $\ord \Jac _f $ are constant on every $B_i$.
Then from Lemma \ref{key lemma 2}, we have
\[
 \int _{A_i} F \LL ^{\fs_\cX} d\mu_\cX = \int _{B_i} (F\circ f_\infty) \LL^{-\ord \Jac _f + \fs_\cY} 
 d\mu_\cY .
\]
Hence
\begin{align*}
&\int _{A} F \LL ^{\fs_\cX} d\mu_\cX \\
&= \sum_{i\in \NN} \int _{A_i} F \LL ^{\fs_\cX} d\mu_\cX \\
& =\sum_{i \in \NN} \int _{B_i} (F\circ f_\infty) \LL^{-\ord \Jac _f + \fs_\cY} 
 d\mu_\cY \\
&= \int _{f_\infty ^{-1}(A)} (F\circ f_\infty) \LL^{-\ord \Jac _f + \fs_\cY} 
 d\mu_\cY .
\end{align*}
\end{proof}

\subsubsection{Key lemmas}

Let $\cX$ be a DM stack of pure dimension $d$.
We define the Jacobian ideal sheaf  $\Jac_\cX$ to be the $d$-th Fitting ideal of
$\Omega _{\cX/k}$. If $\cX$ is smooth, then $\Jac _\cX = \cO_\cX $. 
If $\cX$ is  reduced, then since $k$ is perfect, $\cX$ is generically smooth.
It follows that $\Omega _{\cX/k}$ is generically free of rank $d$ and
the support of $\cO_\cX /\Jac _\cX$ is of positive codimension. 

Let $x \in \cX (K)$. Since for every $ 0 \le n \le \infty$, the natural morphism $\cJ_n \cX \to \cX $
is representable and affine, 
the fiber $( \cJ_n \cX )_x := \cJ_n \cX \times _{\cX,x} \Spec K$ is an affine scheme.

The following lemmas are generalizations of \cite[Lemmas 1.17 and 3.5]{DL-quotient}.
In proving these, we adopt ring-theoretic arguments as in \cite{looi}.

\begin{lem}\label{(n,n-e)-lemma}
Let $f:\cY \to \cX$ be as in Theorem \ref{thm-transformationstackstack}.
Let $\beta_0,\beta_1: D _{\infty,K}   \to \cY$ be (non-twisted)
 arcs over a $K$-point $y \in \cY(K)$, that is, $\beta _0 , \beta _1 \in (J_\infty \cY)_ y(K)= (\cJ_\infty ^1 \cY)_y(K)$.
Assume that $ \pi _n f_\infty (\beta_0) = \pi _n f_\infty (\beta_1)$ for some $n \in \ZZ _{\ge 0}$,
and 
\begin{align*}
e &:= \ord \Jac _ f (\beta_1) < n ,\\
b &:= \ord \Jac _\cX (f_\infty (\beta_1)) < n.
\end{align*}
Then $\pi _{n-  e  } (\beta_0) = \pi _{n-  e}(\beta_1) \in (J_{n-e} \cY)_y (K)$.
\end{lem}

\begin{proof}
Let  $x:= f(y )$.
From the assumption, the homomorphisms 
\[
 (f\beta_0)^* ,  (f\beta_1)^*: \cO _{\cX,x} \to K[[t]]
 \]
  are identical modulo $\fm^{n+1}$. Here $\cO _{\cX,x}$ denotes the complete local ring at $x$ and
 $\fm :=(t)$ is the maximal ideal. 
Hence the homomorphisms $(f\beta_0)^* $ and $  (f\beta_1)^*$
give the same $\cO_{\cX,x}$-module structure to the $K[[t]]$-module $\fm^{n+1}/\fm^{2(n+1)}$.
The  map
\[
 \cO _{\cX,x} \xrightarrow{(f\beta_1)^* -(f\beta_0)^* } \fm^{n+1} \twoheadrightarrow 
 \fm^{n+1}/\fm^{2(n+1)}
\]
is a $k$-derivation, which induces a $K[[t]]$-module homomorphism 
\[
\delta(f\beta_1 ,f\beta_0):(f\beta_1 )^* \Omega _{\cX/k} \to \fm^{n+1}/\fm^{2(n+1)} \twoheadrightarrow 
\fm^{n+1} / \fm^{n+2} .
\]
This annihilates the torsion part of $(f \beta_1 )^* \Omega _{\cX/k}$:
Since the $d$-th Fitting ideal of
$(f\beta_1 )^* \Omega _{\cX/k}$ is $(t^{b})$,  
the torsion part of $(f \beta_1 )^* \Omega _{\cX/k}$ is of length $b$. 
Hence the image of the torsion part in $\fm ^{n+1}/\fm^{2(n+1)}$ is contained in 
 $\fm ^{2(n+1)-b}/\fm^{2(n+1)}$.
Since $n  > b $, we have
\[
2(n+1) - b > n +2.
\]
Thus $\fm ^{2(n+1)-b}/\fm^{2(n+1)}$ is killed by  $\fm^{n+1}/\fm^{2(n+1)} \twoheadrightarrow 
\fm^{n+1} / \fm^{n+2} $.

Consider an exact sequence of $K[[t]]$-modules
\[
 (f\beta_1)^* \Omega_{\cX/k} \to (\beta_1)^* \Omega _{\cY/k} \to (\beta_1)^* \Omega _{\cY/\cX} \to 0 . 
\]
Since $(\beta_1)^* \Omega _{\cY/\cX}$ is of length $e$,  $\delta(f\beta_1 ,f\beta_0)$ lifts to a homomorphism 
\[
 \delta(\beta _1 ,\beta_2) :(\beta_1)^* \Omega _{\cY /k} \to \fm^{n-e+1}/\fm^{n+2},
\]
which derives from
\[
 (\beta_1)^* - (\beta _2) ^* :  \cO _{\cY,y} \to  \fm^{n-e+1}/\fm^{n+2}
 \] 
for some $\beta _2 \in (J_\infty \cY )_y(K)$.
Then
\begin{align*}
\beta_1 &\equiv \beta_2   \mod \fm ^{n-e+1} , \ \text{and} \\
f\beta_0 &\equiv f\beta_2 \mod \fm ^{n+2}.
\end{align*}

Applying the argument above to $\beta _0$ and $\beta _2$,
we obtain $\beta_3$ such that 
\begin{align*}
\beta_2 &\equiv \beta_3   \mod \fm ^{n-e+2}\\
(\text{hence}\ \beta_1 &\equiv \beta _3 \mod \fm ^{ n-e+1 }) , \ \text{and} \\
f\beta_0 &\equiv f\beta_3 \mod \fm ^{n+3}.
\end{align*}
Repeating this, we obtain a sequence $ \beta _i $, $i \in \NN$ such that
\begin{align*}
\beta_1 &\equiv \beta_i   \mod \fm ^{n-e+1} , \ \text{and} \\
f\beta_0 &\equiv f\beta_i \mod \fm ^{n+i}.
\end{align*}
If we put $\beta_\infty$ to be the limit of this sequence, 
 we have 
\begin{align*}
\beta_1 & \equiv \beta _\infty \mod \fm ^{n-e+1} \   \ \text{and} \\
f\beta_0 &= f\beta_\infty .
\end{align*}
Since $f$ is birational and separated, from the valuative criterion \cite[Proposition 7.8]{LMB}, 
$\beta_0$ and $\beta_\infty$ are actually the same.
It follows that $\pi _{n- e} (\beta_0) = \pi _{n- e }(\beta_1) $.
\end{proof}

\begin{lem}\label{keylemma}
Let $f : \cY \to \cX$ be as in Theorem \ref{thm-transformationstackstack}.
Let $K \supset k$ be an algebraically closed field,
$y \in \cY(K)$ and 
$\beta \in (\cJ^l _\infty \cY)_y (K)$.
Suppose that $\ord \Jac _f (\gamma) = e < \infty$, $s_1:= \fs _\cY (\beta)$, 
$s_2:=\fs _\cX (f_\infty \beta)$. 
Assume that $m:=n- \lceil e \rceil >e$ and
$\ord \Jac _\cX (\beta) <  n$. 
Then $f_n^{-1} f_n \pi_n(\beta) \cong \AA_{K}^{e-s_1+s_2}$.
\end{lem}

\begin{proof}
Let $\beta': D_{\infty,K} \to \cY$ be the composite of $\beta$ and the canonical atlas
$D_{\infty , K} \to \cD_{\infty,K}^l$.
Let $m := n - \lceil  e \rceil $ and 
let $\beta _1 \in (\cJ^l _\infty \cY)_y(K)$ be a twisted arc with $ \pi _{m} (\beta )= \pi_{m} (\beta_1)$.
Then as in Lemma \ref{(n,n-e)-lemma}, we obtain a  $\bmu_l$-equivariant derivation
\begin{equation*}
 (\beta')^*-(\beta'_1)^* :  \cO_{\cY,y} \to \fm ^{ml +1 } / \fm ^{nl+1} 
\end{equation*}
and an associated $\bmu_l$-equivariant homomorphism
\begin{equation*}
 \delta (\beta_1) : (\beta') ^* \Omega _{\cY/k}  \to \fm ^{ml +1 } / \fm ^{nl+1} 
\end{equation*}
Let $\beta _2 \in \cJ^l _\infty \cY$ be another twisted arc with $ \pi _{m} (\beta )= \pi_{m} (\beta_2)$.
Then $ \delta (\beta_1) = \delta (\beta_2) $ if and only if $\pi _{n}(\beta_1 ) =\pi _{n}(\beta_2) $.
Therefore if $ \pi_{m}^n $ denotes the natural morphism $\cJ_n \cY \to \cJ_{n '} \cY$, then 
we can regard $(\pi_{m}^n)^{-1}(\pi_{m}(\beta))$ as a subset of
$\Hom _{K[[t]]}^{\bmu_l} ((\beta') ^* \Omega _{\cY/k}  , \fm ^{ml +1 } / \fm ^{nl+1}  )$. 
Here $\Hom _{K[[t]]}^{\bmu_l} (,)$ is the set of $\bmu_l$-equivariant $K[[t]]$-homomorphisms.
In fact, 
the dimensions of $(\pi_{m}^n)^{-1}(\pi_{m}(\beta))$ and $\Hom_{K[[t]]}^{\bmu_l} ((\beta') ^* \Omega _{\cY/k}  , \fm ^{ml +1 } / \fm ^{nl+1}  )$ are equal, and we can identify the two spaces.

Let $x := f(y)$ and $F \subset \Aut (y)$ the largest subgroup
acting on $\cO_{\cY , y} $ trivially. 
Then since $f$ is birational, the natural map $F \hookrightarrow \Aut (y) \to \Aut (x)$ is injective. 
Let  $ b :\bmu_{l} \hookrightarrow \Aut (y)$ and $a : \bmu_{l'} \hookrightarrow \Aut (x)$ be 
embeddings deriving from 
$\beta $ and $f_\infty (\beta)$ respectively. 
We have the following commutative diagram:
\[\xymatrix{
  \zeta \in\bmu_l \ar@{^{(}->}[r] ^b \ar@<2ex>@{>>}[d] \ar@<-2ex>@{|->}[d] & \Aut (y) \ar[d] \\
  \zeta ^{l/l'} \in \bmu_{l'} \ar@{^{(}->}[r]_a & \Aut (x) .
  }
  \]
  Let $\cY'\subset \cY$ be the closed substack where $f$ is not an isomorphism.
From the assumption that $m>e$, $ (\pi_{m}^n)^{-1} \pi _{m}(\beta) \cap (\cJ_ {n} \cY')_y = \emptyset $.
Hence the automorphism groups of any geometric point of $ (\pi_{m}^n)^{-1} \pi _{m}(\beta)$
is the centralizer of $b$ in $F$. It follows that the morphism $ \cJ_n \cY \to \cJ _n \cX $ is
representable around $ (\pi_{m}^n)^{-1} \pi _{m}(\beta)$.
We can regard $f_n^{-1}f_n \pi_n (\beta)$ as a subspace of $(\pi_m^n)^{-1} \pi_m (\beta)\cong \AA^{(n-m)d}_K$.

Consider the natural homomorphism 
\begin{multline}\label{kernel}
 \Hom_{K[[t]]}^{\bmu_l} ( (\beta')^* \Omega _{\cY /k},  \fm ^{ml +1} / \fm^{nl +1} ) \to 
 \Hom_{K[[t]]} ( (f\beta')^* \Omega _{\cX/k},  \fm^{ml +1} / \fm^{nl + 1} ).
\end{multline}
Homomorphisms $\delta (\beta_1)$ and $\delta (\beta_2)$ maps to the same element by this homomorphism
if and only if $f_n \pi_n (\beta _1)=f_n \pi_n (\beta _2) $.
Hence  $f_n^{-1}f_n\pi _n (\beta)$ is isomorphic to the kernel of (\ref{kernel}).

Choose an isomorphism $\cO_{\cY,y} \cong K[[y_1, \dots , y_d]]$ such that 
by the induced $\bmu_l$-action on $K[[y_1, \dots , y_d]]$,
$\zeta \in \bmu_l$ sends $y_i$ to $\zeta ^{b_i} y_i$, $1 \le b_i \le l$.
Then
\[
s_1 = \frac{1}{l} \sum _{i=1}^d (l - b_i) .
\]
We have an isomorphism $(\beta')^* \Omega _{\cY /k} \cong \bigoplus _i K[[t]] dy_i$.
We define a $K[[t^l]]$-homomorphism 
\[
\psi : \bigoplus_i K[[t^l]] dy_i \to (\beta')^* \Omega _{\cY /k} ,
\]
 by  $\psi ( dy_i  ) := t^{l-b_i} dy_i $. The image of this monomorphism
  is the submodule of the $\bmu_l$-invariant
 elements.
 Let $\fn := (t^l)$ be the maximal ideal of $K[[t^l]]$.
Then the  map deriving from $\psi$,
\[
\Hom_{K[[t]]}^{\bmu_l} ((\beta') ^* \Omega _{\cY/k}  , \fm ^{ml +1 } / \fm ^{nl+1} ) \to
\Hom _{K[[t^l]]}( \bigoplus_i K[[t^l]] dy_i , \fn ^{m +1 } / \fn ^{n+1}  ) ,
\] 
is bijective.

If $\cX$ is smooth, we choose an isomorphism
$ \cO_{\cX,x} \cong K[[x_1 ,\dots, x_d]] $
such that by the induced $\bmu _{l'}$-action on  $K[[x_1 ,\dots, x_d]]$, $\zeta \in \bmu_{l'}$ sends
$x_i$ to $\zeta ^{a_i} x_i$, $ 1 \le a_i \le l'$.
Then we have
\[
s_2 = \frac{1}{l'}\sum_i (l'-a_i) = \frac{1}{l}\sum (l - a_i l/l' ).
\]
We have an isomorphism $(f\beta')^* \Omega _{\cX /k} = \bigoplus_i K[[t]]dx_i$.
We define a $K[[t^l]]$-homomorphism
\[
\phi : \bigoplus _i K[[t^l]] dx_i \to (f\beta')^* \Omega _{\cX/k}, 
\]
by $ \phi ( dx_i  ) := t^{ l- a_i l/l'} dx_i $.
Next, suppose that $\cX$ is a variety.
Then $(f\beta')^* \Omega _{\cX /k} \cong (\bigoplus _{i=1} ^d K[[t]] dx_i) \oplus (tors)$ with
$dx_i$ symbols, where
$(tors)$ denotes the torsion part. We define a $K[[t^l]]$-homomorphism
\[
\phi : \bigoplus _i K[[t^l]] dx_i \to  (\bigoplus _{i=1} ^d K[[t]] dx_i) \oplus (tors) \cong (f\beta')^* \Omega _{\cX/k}, 
\]
by $dx_i \mapsto dx_i$.
In both cases, the image of $\phi$ is contained in the submodule of the $\bmu_{l}$-invariant elements.
Here $\bmu_l$ acts on $(f\beta')^* \Omega _{\cX/k}$ through $\bmu_l \twoheadrightarrow \bmu_{l'}$.
Since the natural homomorphism $ (f\beta ')^* \Omega _{\cX/k}  \to (\beta')^* \Omega _{\cY /k}$ is
$\bmu_l$-equivariant, the image of the composite map
\[
 \bigoplus _i K[[t^l]] dx_i \xrightarrow{\phi} (f\beta ')^* \Omega _{\cX/k}  \to (\beta')^* \Omega _{\cY /k}
 \]
is contained in the submodule of $\bmu_l$-invariant elements and the map lifts to
\[
\tau: \bigoplus _i K[[t^l]] dx_i \to \bigoplus _i K[[t^l]] dy_i .
\]

Consider the homomorphism induced by $\tau$,
\[
 \Hom (\bigoplus _i k[[t^l]] dy_i , \fn^{m+1}/\fn^{n+1}) \to 
 \Hom (\bigoplus _i k[[t^l]] dx_i, \fn^{m+1}/\fn^{n+1}).
\]
Its kernel is isomorphic to the kernel of (\ref{kernel}) and to
\[
 \Hom (\Coker \tau , \fn^{m+1}/\fn^{n+1}) .
\]
Hence we have
\begin{align*}
&\dim f_n^{-1}f_n\pi _n (\beta) \\
&= \dim \Hom (\Coker \tau , \fn^{m+1}/\fn^{n+1}) \\
&= \dim \Coker \tau \\
&= \frac{1}{l} ((\dim (\beta ')^*\Omega_{\cY/\cX}) - \sum _i (l- b_i)  + \sum _i (l -a_il/l') )  \\
&= e  -s_1 + s_2 .
\end{align*}
Thus we have proved Lemma \ref{keylemma}.
\end{proof}

\begin{lem}\label{key lemma 2}
Let $B \subset \cJ_\infty \cY$ be a  subset. Suppose that  $\ord \Jac_f |_B \equiv e <\infty$,
$\fs _\cY |_B \equiv s_1$ and $\fs_\cX|_{f_\infty (B)} \equiv s_2$.
Then $B$ is a cylinder if and only if $f_\infty (B)$ is so. If $B$ is a cylinder, then
 we have
\[
 \mu_\cY (B) = \mu_\cX(f_\infty (B)) \LL^{e-s_1 +s_2} .
\]
\end{lem}

\begin{proof}
From Proposition \ref{birational-bijective}, the assumption that $\ord \Jac_f |_B < \infty$ means that
the map $B \to f_\infty (B)$ is bijective.
Since we have assumed that $\cY$ is smooth, $f$ is not an isomorphism all over the singular locus $\cX_{\sing}$ of $\cX$.
Hence $f_\infty (B)$ lies outside $| \cJ _\infty \cX_{\sing}  |$.
In particular, $\ord \Jac _\cX $ takes finite values on $f_\infty (B)$.
Either if $B$ is a cylinder or if $f_\infty (B)$ is cylinder, then $\ord \Jac _\cX |_{f_\infty (B)}$
is bounded from above: In the case where $B$ is a cylinder, consider $(\ord \Jac _\cX)\circ f_\infty =
\ord f^{-1 }(\Jac _\cX)$.

If $f_\infty (B)$ is an $n$-cylinder for $n \gg 0$, then for any point $q \in \pi_n (B)$,
\[
 \pi _n^{-1} (q) \subset \pi_n^{-1} f_n^{-1} f_n (q) \subset f_\infty ^{-1} \pi_n^{-1} f_n(q) 
 \subset f_\infty ^{-1} f_\infty (B) =B.
\]
Since $f_n^{-1}B$ is a constructible subset, it is an $n$-cylinder. 

Next, assume that $B$ is an $(n- \lceil e \rceil )$-cylinder and
that $n$ is large enough to satisfy the condition in Lemma \ref{(n,n-e)-lemma}.
From Lemma  \ref{(n,n-e)-lemma}, for a point $p \in f_n \pi_n (B)$, 
$ \pi_{n-\lceil e \rceil } f_\infty^{-1} \pi_n ^{-1}(p) $ is one point.  Therefore we have
\[
 f_\infty^{-1} \pi_n ^{-1}(p) \subset B 
\]
and 
\[
 \pi_n ^{-1}(p) \subset f_\infty (B). 
\]
From Chevalley's theorem \cite[Th\'eor\`eme 5.9.4]{LMB}, $f_n \pi_n (B)$ is a constructible subset
and hence $f_\infty (B)$ is an $n$-cylinder.

If $B$ is cylinder, from Lemma \ref{keylemma}, we have
\[
 \mu_\cY (B) = \mu_\cX(f_\infty (B)) \LL^{e-s_1 +s_2} .
\]
\end{proof}

\section{Birational geometry of Deligne-Mumford stacks}\label{birationalofstacks}

In this section,  DM stacks are supposed to be reduced and of finite type.

\subsection{Divisors and invariants of pairs}

Let $\cX$ be a DM stack. A {\em prime divisor} on $\cX$ is just a reduced closed substack of
$\cX$ of codimension one. A {\em divisor} (resp.\ {\em $\QQ$-divisor}) is a linear combination of prime divisors
with integer (resp.\ rational) coefficients. 
A divisor is said to be {\em Cartier} if the corresponding divisor on an atlas of $\cX$ is a Cartier divisor.
A divisor $D$ or a $\QQ$-divisor 
is said to be {\em $\QQ$-Cartier} if there exists a positive integer $m$ such that $mD$
is Cartier. If $\cX$ is smooth, a divisor (resp.\ $\QQ$-divisor) is always Cartier (resp.\ $\QQ$-Cartier).
For a Cartier divisor $D$, we can define an invertible sheaf $\cO_\cX(D)$ as follows: For an \'etale morphism
$U \to \cX$ with $U$ scheme, we put $\cO_\cX (D)_U := \cO_U(D_U) $. Here $D_U$ is the pull-back of $D$ to $U$.
The pull-back of a Cartier or $\QQ$-Cartier divisor by a morphism is defined in a obvious way.

Now suppose that $\cX$ is smooth of pure dimension $d$. 
We associate to each $\QQ$-divisor $D$ on $\cX$ a measurable function
$\fI_D : |\cJ_\infty \cX| \to \QQ$ as follows:
If $D$ is a prime divisor and if $\cI_D$ is the defining ideal sheaf of $D$,
then we define 
\[
 \fI_D := \ord \cI_D .
\]
For a general $D$, if we write $D = \sum u_i D_i$ with $D_i$ prime divisor and $u_i \in \QQ \setminus \{0\}$, then we define
\[
 \fI_D := \sum_i u_i \fI_{D_i} .
\]
This function is defined outside a negligible subset $| \cJ_\infty (\bigcup D_i)| $.
For $\QQ$-divisors $D$ and $E$, we have $\fI _{D + E} = \fI _ D + \fI _E$ at least outside a negligible subset.

\begin{defn}
Let $\cX$ be  a  smooth DM stack, $D$ a $\QQ$-divisor of $\cX$ and
$W \subset |\cX|$ a constructible subset. Then we define an invariant
\[
 \Sigma _W( \cX , D) := \int _ {\pi^{-1}(W)} \LL^{\fI_D + \fs _\cX} d \mu _\cX
 \in \bar \fR^{1/r} .
\]
 Here $\pi : |\cJ_\infty \cX| \to | \cX|$ is the natural projection.
\end{defn}

The invariant lies in $\bar \fR^{1/r}$ for a {\em suitable} $r$.
Below, we will take a suitable $r$ each time 
and not mention it hereinafter.

For a smooth DM stack $\cX$, 
the {\em canonical sheaf} $\omega _{\cX}$ is defined to be the sheaf $ \bigwedge ^d \Omega _{\cX/k}$
of differential $d$-forms.
This is an invertible sheaf. 
If $ f : \cY \to \cX$ is a proper birational morphism of smooth DM stacks,
then there is a natural monomorphism $f^* \omega _\cX \to  \omega _\cY$.
There exists an effective divisor $ K _{\cY/\cX}$ on $\cY$ with support in the exceptional locus
such that $\omega_\cY \cong f^* \omega _\cX \otimes \cO_\cY(K_{\cY/\cX})$.
We call $K_{\cY/\cX} $ the {\em relative canonical divisor}.
The defining ideal of $K_{\cY/\cX}$ is nothing but the Jacobian ideal $\Jac _f$ of $f$.

There exists also a {\em canonical divisor} $K_\cX$, that is,
a divisor such that $\cO_\cX (K_\cX) \cong  \omega_\cX$:
Let $X$ be the coarse moduli space of $\cX$ endowed with a morphism $f :\cX \to X$,
$X_0 \subset X$  the smooth locus and $\cX_0 \subset \cX$ the inverse image of $X_0$.
Since $X$ is normal,  $X \setminus X_0$ is of codimension $\ge 2  $.
Let $K_{X_0}$ be a canonical divisor of $X_0$.
Then if we put $K_{\cX_0} := f^* K_{X_0} + K_{\cX_0/X_0}$, then 
$K_{\cX_0}$ is a canonical divisor of $\cX_0$. 
The unique extension of $K_{\cX_0}$ to the whole $\cX$ is a canonical divisor of $\cX$.
Thus a canonical divisor exists.
For a morphism of smooth DM stacks $f : \cY \to \cX$, we have a equation
\[
 K_\cY \equiv f^* K_\cX + K_{\cY/\cX}.
\]

If $\cX$ is not smooth but only normal and $\cX_0 \subset \cX$ be the smooth locus, then
a canonical divisor $K_\cX$ of $\cX$ is defined to be a divisor such that 
$K_{\cX}|_{\cX_0}$ is a canonical divisor of $\cX_0$.

\begin{thm}\label{MainTheorem}
Let
\[\xymatrix{
& \cY \ar[dl]_f \ar[dr]^{f'} &\\
\cX&&\cX' 
} 
\]
be a diagram
consisting of smooth DM stacks of pure dimension and tame proper birational morphisms.
Let $D$ and $D'$ be $\QQ$-divisors on $\cX$ and $\cX'$ respectively and
let $W$ and $W'$ be constructible subsets of $|\cX|$ and $|\cX'|$ respectively.
Then if $f^{-1}(W) = (f')^{-1}(W')$ and 
\[
f^* D - K_{\cY/\cX} = (f')^* D' - K_{\cY/\cX'},
\]
then
\[
 \Sigma _W(\cX ,D)= \Sigma _{W'}(\cX',D')  .
\]
\end{thm}

\begin{proof}
Since
\[
 \fI_D \circ f_\infty - \ord \Jac_f = \fI_{D'} \circ f'_\infty - \ord \Jac_{f'} ,
\]
 the assertion follows from Theorem \ref{thm-transformationstackstack}.
\end{proof}

\begin{cor}\label{corollaryRuan1}
With the notation of Theorem \ref{MainTheorem}, if $K_{\cY/\cX} = K _{\cY/\cX'}$
 and $f^{-1}(W) = (f')^{-1}(W')$, then we have 
\[
    \sum_{\cV \subset \cJ_0 \cX} \{ |\cV|_W \} \LL ^{\sht (\cV)} 
    = \sum _{\cV' \subset \cJ_0 \cX'} \{ |\cV'|_{W'} \} \LL ^{\sht (\cV')}.
\]
Here the sums run over the connected components and $|\cV|_W$ denotes the inverse image of $W$ in $|\cV|$.
\end{cor}

\begin{proof}
It follows from the fact that
\[
 \Sigma _W(\cX,0)=\sum_{\cV \subset \cJ_0 \cX} \{ |\cV|_W \} \LL ^{\sht (\cV)} .
\]
\end{proof}

\subsection{Homological McKay correspondence and discrepancies}

Suppose $k=\CC$. Let $X$ be a $\QQ$-Gorenstein variety, that is, a normal variety
with $\QQ$-Cartier canonical divisor $K_X$.  
Let $f:Y \to X$ be a resolution of $X$, that is, $Y$ is smooth and $f$ is proper and birational.
Then we can attach a rational number $a(E,X)$ to each exceptional divisor $E$ 
such that we have a numerical equivalence
\[
 K_Y \equiv f^* K_X + \sum _{E:\text{exceptional}} a(E,X) E.
\]
We call $a(E,X)$ the {\em discrepancy} of $E$ with respect to $X$.
If $a(E,X)$ is zero for every $E \subset Y$, then $Y$ is said to be a {\em crepant resolution}
of $X$. 
The {\em discrepancy} of $X$, denoted $\discrep (X)$, is defined as follows:
\[
 \discrep (X) := \inf \{a(E,X) |  Y \to X, E \subset Y \ \text{exceptional}\}.
\]  
Here $Y \to X$ runs over all resolutions of $X$.
It is well-known that either $\discrep (X) = - \infty$ or $\discrep (X) \ge -1$.
If $\discrep (X) \ge -1$, then $\discrep (X) $ is equal to the minimum of 
$a(E,X)$ for exceptional divisors $E$ on a {\em single} resolution such that the exceptional locus is simple
normal crossing.
The discrepancy is an important invariant of singularities in the minimal model program.

Consider a finite subgroup $G$ of $GL_d(\CC)$ and the quotient variety $X := \CC^d/G$.
Then $X$ is $\QQ$-Gorenstein and has log terminal singularities, that is, $\discrep (X) > -1$. Moreover if $G \subset SL_d(\CC)$, then $X$ is Gorenstein and has 
canonical singularities, that is, $\discrep (X) \ge 0$.
Let $g \in G$ be an element of order $l$ and let $\zeta_l = \exp (2\pi \sqrt{-1}/l)$.
Choosing suitable basis of $\CC^d$, we write
\[
 g = \diag(\zeta_l ^{a_1} , \dots , \zeta _l ^{a_d} ) ,\ 0 \le a_i \le l-1.
\] 
Then we define the {\em age} of $g$ to be 
\[
\age (g) := \frac{1}{l}\sum_{i=1}^d a_i \in \QQ. 
\]
If $g \in \SL_d(\CC)$, then $\age (g) $ is an integer.
Now we can deduce, from Corollary \ref{corollaryRuan1}, 
the homological McKay correspondence. 
It was proved by Y.\ Ito and Reid \cite{ito-reid}  for dimension 3 and by Batyrev \cite{non-archi} for arbitrary dimension. See \cite{reid-bourbaki} for a nice survey of this subject.

\begin{cor}[Homological McKay correspondence]\label{HomologicalMcKay}
Suppose that $G \subset \SL _d (\CC)$ and there exists a crepant resolution $Y \to X$. 
For an even number $i \ge 0$,
let $n_i := \sharp \{g \in \Conj (G) | \age(g) = i /2 \}$.
Then 
\[
 H ^{i} (Y,\QQ ) \cong 
 \begin{cases}
\QQ(-i/2)  ^{\oplus n_i} & (i:\text{even}) \\
0 & (i:\text{odd}) . 
\end{cases}
\]
\end{cor}

\begin{proof}
Let $\cX := [\CC^d/G]$. Then the natural morphism
$\cX \to X$ is a proper birational morphism. Furthermore, since $G\subset \SL_d (\CC)$,
$\cX$ and $X$ are isomorphic in codimension one.  In particular, $K _\cX$ is  the pull-back of $K _X $.
(Hence $\cX$ is a crepant resolution in a generalized sense.)

There exists a smooth DM stack $\cY$ and proper birational morphisms $f:\cY \to \cX$ and $f':\cY \to Y$:
For example, we can take a resolution of the irreducible component of $Y \times_X \cX$ dominating
$Y$ and $\cX$. For resolutions of DM stacks, see Subsection \ref{subsec-singular-stack}.
Then we have $K_{\cY/\cX} =K_{\cY/Y}$.
For $ 0 \le a_i \le l-1$, if we put 
\[
 a_i':= 
 \begin{cases}
l & (a_i =0) \\
a_i & (a_i \neq 0),
\end{cases}
\]
then 
\[
 \frac{1}{l} \sum _i (l-a_i') = d  - \sharp \{ i | a_i =0\}-  \frac{1}{l}\sum _i a_i.
\]
From Corollary \ref{corollaryRuan1}, 
\[
 \{Y\} = \sum _{\cV \subset I\cX} \{ \cV \} \LL^{\sht (\cV)} =
  \sum _{g \in  \Conj (G)} \{[(\CC^d)^g/C_g]  \} \LL ^  {d-\dim (\CC^d)^g - \age (g)} .
\]
Sending these elements to $\hat K_0(MHS)$, we have 
\[
 \chi_h (Y) = \sum_{g \in  \Conj (G)} \LL^{d - \age (g)} .
\]
Since $K_0(MHS) \to \hat K_0 (MHS) $ is injective, the equation holds also in  $K_0(MHS)$.
 Batyrev \cite{non-archi} proved that $H^i _c (Y,\QQ)$ has pure Hodge structure of weight $i$.
Let $m_i:= \{ g \in \Conj (G)| -d + \age (g) = -i/2\}$.
Looking at the weight $ i $ part, we have 
\[
 H _c ^i(Y,\QQ ) \cong 
 \begin{cases}
 0 & (i: \text{odd}) \\
 \QQ ( - i / 2 )^{\oplus m_i} & (i : \text{even}) .
\end{cases}
\]
Here we have used the semisimplicity of polarizable pure Hodge structure. 
The corollary follows from the Poincar\'e duality.
\end{proof}

Next, for general $G \subset GL _d (G)$ without reflection, we deduce an expression of the discrepancy of $X= \CC^d/G$ in terms of
ages of $g \in G$. 

\begin{cor}
Suppose that $G$ contains no reflection.
We have an equation
\[
 \discrep (X) = \min \{ \age (g) | 1 \ne g \in G\} -1  .
 \]
\end{cor}

\begin{proof}
Let $\cX := [\CC^d/G]$. 
and  $V \subset  |\cX|$ the locus of points with nontrivial automorphism group.
We take a resolution $f :Y \to X$ which is an isomorphism over the smooth locus of $X$.
Suppose that the exceptional locus $W \subset Y$ is simple normal crossing. 
From Corollary \ref{corollaryRuan1}, as in the proof of Corollary  \ref{HomologicalMcKay}, we have
\[
 \Sigma _W (Y , - K_{Y/X}) = \Sigma_V  (\cX , 0) .
\]

The left hand side can be computed explicitly as follows.
Write $  K_{ Y/X} = \sum_ {i \in I}  e_i E_i$. 
Note that $e_i >-1$.
For $ s= (s _ i)_{i \in I} \in (\ZZ_{\ge 0})^I $, 
let $I_s:= \{ i \in I | s_i >0\}$.  For any subset $J \subset I$, we define
$E_{J}^{\circ}:= \bigcap_{i \in J} E_i \setminus \bigcup _{i \in I \setminus J} E_i$. 
Then we have
\[
 \mu _Y ( \bigcap _{i \in I} \cI _{D_i}^{-1} ( s_i) ) = 
 \{ E_ {I_s} ^\circ \} \LL^{- \sum s_i} (\LL-1)^{\sharp I_s}.
\]
(See \cite[The proof of Theorem 2.15]{Craw-introduction-motivic}. Note that our definition of the motivic measure differs from that of \cite{Craw-introduction-motivic} by the multiplication of $\LL^d$. See also
Remark \ref{rem-difference-Ld}.)
Hence
\[
 \Sigma _W (Y , - K_{Y/X}) = \sum _{0 \ne s \in  (\ZZ_{\ge 0})^I }  \{ E_ {I_s} ^\circ \} 
 \LL^{- \sum (e_i +1) s_i} (\LL-1)^{\sharp I_s} . 
\]
Therefore the dimension of the right hand side is equal to
\begin{align*}
 &\dim ( \sum _{i \in I} \{ E_ {\{i\}} ^\circ \} \LL^{ - e_i -1} (\LL-1) ) \\
&= d - 1 + \max \{ -e_i  | i \in I\} \\
&= d -1 - \discrep (X) . 
\end{align*}
  On the other hand, the dimension of $\Sigma _V (\cX,0)$
is equal to $ d - \min \{\age (g) | 1 \ne g \in G\} $.
This proves the assertion.
\end{proof}

\subsection{Orbifold cohomology}

Chen and Ruan \cite{CR} constructed a new kind of cohomology, called the orbifold cohomology, 
for topological orbifolds. 
We define the orbifold cohomology for DM stacks, the algebraic counterpart of the topological orbifold.

\begin{defn}
\begin{enumerate}
\item Let $\cX$ be a smooth DM stack over $\CC$. We define the orbifold cohomology of $\cX$ as follows:
For each $i \in \QQ$, 
\[
 H^i  _{orb}(\cX ,  \QQ) := \bigoplus _{\cV \subset \cJ_0 \cX} H^{i -2 \sht (\cV)} (\bar \cV, \QQ)
 \otimes \QQ (- \sht (\cV)) .
\]
Here by convention, we put $H ^i (X,\QQ ) = 0 $ for a variety $X$ and $ i \notin \ZZ$.

\item Let $\cX$ be a smooth DM stack over a finite field $k$.
Let $r$ be the least common multiple of $\sht (\cV)$, $\cV \subset \cJ_0 \cX$.
Suppose that $\QQ_p(1/r)$ exists. (This holds after replacing $k$ with its finite extension.)
We define the $p$-adic orbifold cohomology of $\cX \otimes \bar k$ as follows:
For each $i \in \QQ$, 
\[
 H^i  _{orb}(\cX \otimes \bar k, \QQ_p) := \bigoplus _{\cV \subset \cJ_0 \cX} H^{i -2 \sht (\cV)} (\bar \cV \otimes \bar k,\QQ_p)
 \otimes \QQ_p (- \sht (\cV)) .
\]
\end{enumerate}

\end{defn}

\begin{lem}
Let $\cX$ be a proper smooth DM stack over a perfect field $k$.
\begin{enumerate}
\item If $k= \CC$, then $H^i _{orb}(\cX ,\QQ)$ is a pure Hodge structure of weight $i$.
\item If $k$ is a finite field and if for $r \in \NN$ as above, $\QQ_p(1/r)$ exists, then 
$H^i _{orb}(\cX \otimes \bar k ,\QQ_p)$ is pure of weight $i$.
\end{enumerate}
\end{lem}

\begin{proof}
1. Each connected component $\cV$ of $\cJ_0 \cX$ is a proper smooth DM stack.
Then the coarse moduli space $\bar \cV$ is a proper variety with quotient singularities. 
Therefore $H^i (\bar \cV , \QQ)$ is a pure Hodge structure of weight $i$ and 
$H^i _{orb} (\cX , \QQ)$ is also so.

2. 
We claim that the constant sheaf $\QQ_p$ on $\bar \cV \otimes \bar k$ is pure of weight zero. Since the purity of sheaves is a local property, 
we may assume $\bar \cV = M/G$ for some smooth variety $M$
and a finite group $G$ acting on $M$. Then since $q :M \to M/G$ is finite, $q_* \QQ_p$ is pure of weight zero \cite{Del-Weil-II} and $\QQ_p = (q_* \QQ_p)^G$ is also so.

Since $\bar \cV$ is proper, $H^i(\bar \cV \otimes \bar k, \QQ_p)$ is pure of weight $i$
and $H_{orb}^i (\cX \otimes \bar k , \QQ_p)$ is also so.
\end{proof}

The following corollary was conjectured by Ruan \cite{ruan}.
 A weak version was proved by Lupercio and Poddar \cite{Lupercio-Poddar}
 and the author \cite{Twistedjet} independently.

\begin{cor}\label{ruan'sconjecture3}
Let $\cX$ and $\cX'$ be proper smooth DM stacks over $k=\CC$.
Assume that there exist a smooth DM stack $\cY$ and proper birational morphisms $f:\cY \to \cX$ and 
$f': \cY \to \cX$ such that $K_{\cY/\cX} = K_{\cY /\cX'}$.
Then we have an isomorphism of Hodge structures
\[
H^i  _{orb}(\cX , \QQ) \cong H^i  _{orb}(\cX' , \QQ), \ \forall i \in \QQ .
\]
(We do not assert that there exists a natural isomorphism.)
\end{cor}

\begin{proof}
From Corollary \ref{corollaryRuan1}, we have
\begin{equation}\label{equation-orbifold}
  \sum_{\cV \subset \cJ_0 \cX} \chi_h ( \bar \cV ) \LL ^{\sht (\cV)} 
    = \sum _{\cV' \subset \cJ_0 \cX'}\chi_h (\bar  \cV')  \LL ^{\sht (\cV')}.
\end{equation}
Define 
\[
 H^{ev,i}  _{orb}(\cX ,  \QQ) := \bigoplus _{\substack{\cV \subset \cJ_0 \cX \\ i -2 \sht (\cV) : \text{even}}} H^{i -2 \sht (\cV)} (\bar \cV, \QQ)
 \otimes \QQ (- \sht (\cV)) ,
\]
and $H^{odd,i}_{orb}(\cX , \QQ)$ similarly. Then 
\[
H^i_{orb}(\cX,\QQ) =H^{ev,i}_{orb}(\cX,\QQ)\oplus H^{odd,i}_{orb}(\cX,\QQ) .
\]
Looking at the weight $i$ part of (\ref{equation-orbifold}), we have
\[
 [H^{ev,i}_{orb}(\cX,\QQ)] - [ H^{odd,i}_{orb}(\cX,\QQ) ] 
 =[H^{ev,i}_{orb}(\cX',\QQ)] - [ H^{odd,i}_{orb}(\cX',\QQ) ]. 
\]
Since the two terms of each hand side do not cancel out, we have
$ [H^{ev,i}_{orb}(\cX,\QQ)]=[H^{ev,i}_{orb}(\cX',\QQ)]$ and $[H^{odd,i}_{orb}(\cX,\QQ)]=[H^{odd,i}_{orb}(\cX',\QQ)] $. 
Hence $[H^i_{orb}(\cX,\QQ)] = [H^i_{orb}(\cX',\QQ)]$.

We claim that for an arbitrary proper smooth  DM stack $\cX$ over $\CC$, 
$H^{i}_{orb}(\cX,\QQ)$ is semisimple.
A polarization of 
$H^{i-\sht (\cV )} (\bar \cV, \QQ)$ induces a non-degenerate bilinear form $Q_{\bar \cV}$ on
$H^{i-\sht (\cV)} (\cV , \QQ ) \otimes \QQ(-\sht(\cV))$ for which
the $\frac{1}{r}\ZZ$-indexed Hodge decomposition of $ (H^{i-\sht (\cV)} (\cV , \QQ ) \otimes \QQ(-\sht(\cV)))$
is orthogonal. 
We define a bilinear form $Q_{\cX}$ on $H^{i}_{orb}(\cX,\QQ)$
to be  the direct sum of $Q_{\bar \cV}$.
Then $Q_{\cX}$ is non-degenerate and the Hodge decomposition of $H^{i}_{orb}(\cX,\QQ)$
is orthogonal for this.
We can see that $H^{i}_{orb}(\cX,\QQ)$ is semisimple like the usual polarizable
Hodge structure.

Now the equation above, $[H^i_{orb}(\cX,\QQ)] = [H^i_{orb}(\cX',\QQ)]$, implies that
$ H^i_{orb}(\cX,\QQ)\cong H^i_{orb}(\cX',\QQ)$.
\end{proof}

\begin{cor}\label{corollaryRuan2}
Assume that $k$ is a finite field.
Let $\cX$ and $\cX'$ be proper smooth DM stacks
whose $p$-adic orbifold cohomology groups can be defined.
Assume that there exist a smooth DM stack $\cY$ and
tame proper birational morphisms $f:\cY \to \cX$ and 
$f': \cY \to \cX$ such that $K_{\cY/\cX }=K_{\cY/\cX'}$.
Then we have the following 
isomorphisms of Galois representations:
\[
H^i  _{orb}(\cX \otimes _k \bar k , \QQ_p)^{ss} \cong H^i  _{orb}(\cX' \otimes _k \bar k, \QQ_p)^{ss}, \ \forall i \in \QQ .
\]
\end{cor}

\begin{proof}
Similar as Corollary \ref{ruan'sconjecture3} except for the semisimplicity.
\end{proof}

\subsection{Convergence and normal crossing divisors}

\begin{defn}
Let $\cX$ be a smooth DM stack and $D = \sum _{i=1}^n D_i$ a divisor of $\cX$
 with $D_i$ distinct prime divisors.
We say that $D $ is \textit{normal crossing} if 
the pull-back of $D$ to an atlas of $\cX$ is (analytically) normal crossing.
\end{defn}

Let $\cX$ be a smooth DM stack and $x \in \cX(\bar k)$.
Let $\hat \cX = [\Spec \bar 
k [[x_1 ,\dots, x_d]]/G]$ be the completion of $\cX$ at $x$.
If $D$ is a normal crossing divisor on $\cX$,
then for  suitable local coordinates $x_1 , \dots , x_d$,
 the pull-back of $D$ to $\hat \cX$ is defined by 
a monomial $x_1 x_2 \cdots x_c$, $c \le d$.

\begin{defn}
A normal crossing divisor $D$ is said to be \textit{stable normal crossing} if 
for every $x \in \cX(\bar k)$, 
every irreducible component of its pull-back to $ \Spec \bar k [[x_1 ,\dots, x_d]] $  is stable under the $G$-action.
\end{defn}

If $D$ is stable normal crossing and $l \in \NN$ is prime to the characteristic of $k$,
then for each embedding $a : \bmu _l  \hookrightarrow G$,
 we can choose local coordinates $x_1 , \dots , x_d$ so that 
the $\bmu_l$-action is linear and diagonal and th pull-back of $D$ is defined by a monomial simultaneously: 
Suppose that  the pull-back of $D$ is defined by $x_1 \cdots x_ c$. 
Since each irreducible component is stable under the $\bmu_l$-action,
for $1 \le i \le c$ and $\zeta \in \bmu_l$, $\zeta(y_i)$ lies in the ideal $(y_i)$.
Let $\bar \zeta$ be the linear part of $\zeta$, namely
\[
 \zeta (y_i) = \bar \zeta (y_i)+ \text{(terms of order $\ge$ 2)} .
\]
If we replace $y_i$ with
\[
 y_i ' := \sum _{\zeta \in \bmu_l} \bar \zeta ^{-1} \zeta (y_i), 
\]
then the $\bmu_l$-action is linear. We have the identity of ideals $(y_i)=(y_i')$, $1 \le i \le c$,
hence the pull-back of $D$ is still defined by a monomial.
Then the $\bmu_l$-action must be diagonal about coordinates $x_1 ,\dots, x_c$.
Replacing the rest coordinates, we can diagonalize the action.

\begin{rem}
The stable normal crossing divisor is a notion for stacks corresponding the $G$-normal pair in \cite{non-archi}.
\end{rem}

\begin{prop}\label{prop-convergence}
Let $\cX$ be a smooth DM stack, $D = \sum_{i=1}^m u_i D_i$ a $\QQ$-divisor with stable normal crossing support and
$W \subset |\cX|$ a constructible subset.
Then  $ \Sigma_W(\cX,D) \ne \infty $ if and only if
 $u _i < 1$ for every $i$ with $D_i \cap W \ne \emptyset$.
\end{prop}

\begin{proof}
From the semicontinuity of dimension of fibers, 
it suffices to show the proposition in the case where
$W = \{x\}$ with $x \in \cX (\bar k)$. 
Shrinking $\cX$, we can assume that every $D_i$ contains $x$.
Take the completion 
\[
 \hat \cX := [\Spec \bar k[[x_1, \dots, x_d ]]/G]
\]
 of $\cX$ at $ x$. Take an embedding $a : \bmu_l \hookrightarrow G$ and
 choose local coordinates $x_1, \dots , x_d$  so that $\bmu_l$ acts linearly and diagonally and the pull-back of $D_i$ is defined by
 $x_i =0$ for $1 \le  i \le c \le d$ and $x \notin D_i$ for $i > c$. Suppose that
 $\zeta \in \bmu_l$ acts by  $\diag (\zeta ^{a_1},\dots,\zeta^{a_d})$, $1 \le a_i \le l$.

Let $v$ be the $\bar k$-point of $\cJ_0 \cX$ corresponding $(x,a)$.
Then for $0 \le n \le \infty$, 
the fiber  $(\cJ_n \cX)_v$ of $\cJ_n \cX \to \cJ_0 \cX$ over $v$ is identified with 
\[
  \Hom _{\bar k[[t]]}^{\bmu_l} (\bar k[[x_1 ,\dots , x_d]],\bar k[[t]]/t^{nl+1}) .
\]

If $\sigma \in (\cJ_\infty \cX)_v $ and $\sigma (y_i) = \Sigma _{j \ge 0} \sigma _{ij} t^{lj + a_i}$, then
the order of the ideal $(x_i)$ along the twisted arc  $\sigma $ is 
\[
  \frac{a_i}{l} + \min \{ j | \sigma _{ij} \ne 0 \} .
\] 
For a multi-index $s= (s_1 ,\dots s_c) \in (\ZZ_{\ge 0})^c $, we define $V_s \subset (\cJ_\infty \cX)_v$ to be the set 
of $\sigma$ such that $\min \{j | \sigma_{ij}  \neq 0\} = s_i$ for every $1 \le i \le c$.
Then 
\begin{equation}\label{integrationexplicit}
 \int _{\pi_0^{-1} (\bar v)} \LL^{\fI_D +\fs_\cX}d \mu_\cX =\sum _s \mu _{\cY} (\bar V_s) \LL ^{\sum _{i=1}^c u_i (s_i + a_i/l) } .
\end{equation}
Here $\bar v \in |\cJ_0 \cX|$ is the image of $v$ 
and $\bar V_s \subset |\cJ_\infty \cX|$ is the image of $V_s$.
For sufficiently large $n$, 
\[
 \dim \pi_n (\bar V_s) = (d-c)n + \sum_{i=1}^c (n-s_i) = d n - \sum_{i=1}^c s_i.
 \]
Hence 
\[
\dim \mu _{\cY} (\bar V_s) \LL ^{\sum _{i=1}^c u_i (s_i + a_i/l) } =  (u_i-1)s_i +C .
\]
Here $C$ is a constant independent of $s$.
Hence the right hand side of (\ref{integrationexplicit}) converges 
if and only if $u_i < 1$ for every $i$.
\end{proof}

\begin{prop}\label{StablyNormalCrossing}
Let $\cX$ be a smooth DM stack and $D$ a normal crossing divisor of $\cX$.
Then there exist a smooth DM stack $\cY$ and a representable proper birational morphism $f:\cY \to \cX$ 
such that $f$ is 
an isomorphism over  $\cX \setminus D$ and 
 $f^{-1}( D)$ is stable normal crossing.
\end{prop}

\begin{proof}
Let $\Delta$ be an irreducible component of $D$. Note that $\Delta$ can be singular.
For each $\bar k$-point $p \in \Delta$, we define the {\em index} $i_{\Delta}(p)$ as follows:
Let $[ \Spec \bar k[[x_1 , \dots , x_d]] / G]$ be the completion of $\cX$ at $p$.
Choose coordinates so that the pull-back of $\Delta$ is defined by $x_1 \cdots x_c$.
Then $i_{\Delta}(p):=c$.

The function
\[
 i _{\Delta}: \Delta (\bar k) \to \NN
\]
is upper semi-continuous. The locus $V \subset \Delta$ of the points of the maximum index is a closed smooth substack
 defined over $k$.
Let $\cX' \to \cX$ be the blow-up along $V$, $E$ the exceptional divisor, $D'$ the strict transform of $D$
and $\Delta '$ the strict transform of $\Delta$. Then $D' \cup E$ is normal crossing.
For every irreducible component $E_1$ of $E$ and for every $p \in E_1 (\bar k)$,
we have  $i_{E_1} (p)=1$. Moreover the maximum value of the function $i _{\Delta'}$ is less than
that of $i_{\Delta}$. Hence repeating blow-ups, we obtain a proper birational morphism 
$f:\cY \to \cX$ such that every $\bar k$-point of $f^{-1}(D)$ 
 has index one with respect to every irreducible component of $f^{-1}(D)$.
It means that $f^{-1}(D)$ is stable normal crossing. 
\end{proof}

\begin{prop}\label{prop-convergence2}
Let $\cX$ be a smooth DM stack, $D = \sum u_i D_i$ a $\QQ$-divisor with (not necessarily stable) normal crossing support
and $W \subset |\cX|$ a constructible subset. 
Then 
$ \Sigma_W(\cX,D) \ne \infty$
if and only if $u _i < 1$ for every $i$ with $D_i \cap W \ne \emptyset$. 
\end{prop}

\begin{proof}
Let $f:\cY \to \cX$ be a morphism as in Proposition \ref{StablyNormalCrossing}.
Then by Theorem \ref{MainTheorem}, we have
\[
 \Sigma _W(\cX,D)= \Sigma _{f^{-1}(W)}(\cY , f^*D - K_{\cY/\cX}) .
\]
Shrinking $\cX$, we may assume that every $D_i$ meets $W$.
By the standard calculation in the minimal model program, if
$u _i < 1$ for every $i$, then coefficients in $f^*D - K_{\cY/\cX}$ are also all $<1$.
Therefore from Proposition \ref{prop-convergence}, 
the $\Sigma _W(\cX,D)= \Sigma _{f^{-1}(W)}(\cY , f^*D - K_{\cY/\cX}) \ne \infty$.
If $u_i \ge 1$ for some $i$, the coefficient of the strict transform of $D_i$ in $f^*D - K_{\cY/\cX}$
is also $u_i \ge 1$. Hence again from Proposition \ref{prop-convergence}, 
$ \Sigma_W(\cX,D) = \infty$.
\end{proof}

\subsection{Generalization to singular stacks}\label{subsec-singular-stack}

From now on, we assume that the base field $k$ is of characteristic zero.

Thanks to Hironaka \cite{Hironaka}, for every  reduced variety $X$, 
there exists a resolution of singularities, that is, a proper birational morphism 
$Y \to X$ with $Y$ smooth.  
Villamayor \cite{VillamayorConstructiveness} \cite{VillamayorPatsching} constructed
 a resolution algorithm commuting with smooth morphisms.
 See also \cite{BierstoneMilmanCanonicalResolution}, \cite{EncinasVillamayorActa}.
Let $\cX$ be a reduced DM stack and $M$ an atlas. Then we
obtain a groupoid space $N := M \times _\cX M \rightrightarrows M $.
The associated stack of this groupoid space is canonically isomorphic to $\cX$.
Let $\tilde  N $ and $\tilde M$ be smooth varieties obtained from $N$ and $M$ respectively
 by a resolution algorithm commuting with \'etale
morphisms. Then we obtain a groupoid space $\tilde N \rightrightarrows \tilde M$.
Its associated DM stack $\tilde \cX$ is smooth and the natural morphism $\tilde \cX  \to \cX$ is
representable,
proper and birational. 
Thus for every reduced DM stack $\cX$, there exists a representable proper birational morphism
$\cY \to \cX$ with $\cY$ smooth.

Let $\cX$ be a normal DM stack and $D$ a $\QQ$-divisor on $\cX$.
Suppose that $K_\cX + D$ is $\QQ$-Cartier. Then we say that
 the pair $(\cX , D)$ is a {\em log DM stack}. 
For a resolution $f : \cY \to \cX$, we define a $\QQ$-divisor $E$ on $\cY$ by
\[
K_\cY + E \equiv f^* ( K_\cX +D) .
\]

\begin{defn}
Let the notations be as above. Let $W \subset |\cX|$ be a constructible subset. 
Then we define an invariant
\[
\Sigma _W(\cX , D) := \Sigma _{f^{-1} (W)}( \cY ,E)  .
\]
\end{defn}

This is independent of the choice of resolutions, thanks to Theorem \ref{MainTheorem}.

\begin{defn}\label{def-KLT}
We say that the pair $(\cX,D)$ is \textit{Kawamata log terminal (KLT for short)} if
for every representable resolution $\cY \to \cX$, every coefficient of the divisor $E$ defined as above
 is less than one.
\end{defn}

In fact, we can see if $(\cX,D)$ is KLT by looking at only one
resolution with $E$ normal crossing.
The pair $(\cX,D)$ is KLT if and only if 
for an atlas $M \to \cX$ and the pull-back $D'$ of $D$ to $M$, the pair $(M,D')$ is KLT.

\begin{prop}
The invariant $\Sigma _W(\cX , D)$ is not infinite if and only if
$(\cX , D)$ is KLT around $W$ (that is, for some open substack 
$\cX_0 \subset \cX$ containing $W$, $(\cX_0 , D|_{\cX_0})$ is KLT).
\end{prop}

\begin{proof}
It is a direct consequence of Proposition \ref{prop-convergence2} and the definition of the invariant.
\end{proof}

\begin{thm}\label{thmKequiv}
Let $(\cX,D)$ and $(\cX',D')$ be log DM stacks. 
Let $W$ and $W'$ be  constructible subsets of $|\cX|$ and $|\cX'|$ respectively.
Assume that there exist a smooth DM stack $\cY$ and proper birational morphisms $f:\cY \to \cX$ and 
$f': \cY \to \cX$ such that $f^* ( K_\cX + D ) \equiv  f^* ( K_{\cX'} + D' )$
 and $f^{-1}(W) = (f')^{-1}(W')$.
Then we have 
\[
 \Sigma _W(\cX ,D)= \Sigma _{W'}(\cX',D')  .
\]
\end{thm}

\begin{proof}
It follows from  Theorem \ref{MainTheorem}.
\end{proof}

\subsection{Invariants for varieties}
If $X$ is a variety, then we can describe the invariant $\Sigma _W(X,D)$ with 
the motivic integration over $X$ itself.
It  gives us a canonical expression of the invariant.

Let $X$ be a normal variety and $D$ a $\QQ$-divisor on $X$ such that
$m (K_X +D)$ is Cartier  and $mD$ is an integral divisor for some $m \in \NN$. 
Let $\cK$ be the sheaf of total quotient rings of $\cO_X$.
There exists a fractional ideal sheaf $\cG \subset \cK$ such that 
$\cG = \cO_X(-mD)$ outside the singular locus of $X$,
and  
the canonical isomorphism $(\Omega^d_X)^{\otimes m} \cong  \cG \cO_X(m(K_X +D)) $ outside the singular locus
 of $X$ 
 extends
to an epimorphism $ (\Omega _X ^d)^ {\otimes m} \twoheadrightarrow \cG \cdot \cO_X(m(K_X +D)) $ all over $X$.

\begin{defn}
Let $\cI \subset \cK$ be a fractional ideal sheaf.
We define the {\em order function} of $\cI$ as follows;
\begin{align*}
\ord \cI : J_\infty X & \to \ZZ  \cup \{\infty\} \\
(\gamma: \Spec K[[t]] \to X )&\mapsto 
\begin{cases}
n & (\gamma^{-1} \cI = (t^n) \subset K((t))) \\
\infty & (\gamma^{-1} \cI = (0)) .
\end{cases}  
\end{align*}
\end{defn}

If $A$ is a Cartier divisor and $\cI_A= \cO_\cX (-A)$ is the corresponding fractional ideal sheaf, then
we have  $\ord \cI_A = \fI_A$, at least outside a negligible subset.

We have the following expression of the invariant $\Sigma _W(X,D)$.

\begin{prop}
Let $W \subset X$ be a constructible subset.
We have 
\[
\Sigma _W( X ,D) = \int _{\pi _0 ^{-1 } (W)} \LL ^{(1/m)\ord \cG} d\mu_X .
\]
\end{prop}

\begin{proof}
Let $f :Y \to X$ be a resolution of singularities.
We define a $\QQ$-divisor $E$ on $Y$ by 
\[
 K_Y + E \equiv f^* (K_X + D).
\]
 Let $\cI  = \cO _X (-mE)$.
We have natural morphisms
\begin{gather*}
  f^* \cO_X (m(K_X +D)) \cong \cI ^{-1} \omega_Y ^{\otimes m} \\
  f^* \Omega _X ^d \twoheadrightarrow \Jac_f  \omega _Y \\
 (\Omega _X ^d)^ {\otimes m} \twoheadrightarrow \cG \cO_X (m(K_X +D)).
\end{gather*}
Therefore we have
\[
(\Jac _f)^m = f^{-1 } \cG \cdot \cI ^{-1} .
\]
We obtain an equation of measurable functions
\[
  \frac{1}{m} \ord \cG \circ f_\infty - \ord \Jac _f =  \fI _E  .
\]
Now the proposition follows from Theorem \ref{TransformationDenefLoeser}.
\end{proof}

\begin{rem}
This kind of  integration over a singular variety
 was considered in \cite{DL-quotient} (see also \cite{looi})
in relation to the McKay correspondence and 
 in \cite{DimJet} and  \cite{EinMustataYasuda} in relation to the discrepancy of singularities.
\end{rem}

%%%%%%%%%%%%%%%%%%%%%%%%%%%%%%%%%%%%%%%%%%%%%%%%%%%%%%%%%%%%%%%%%%%%%%%%%%%%%%%

\end{document}